\documentclass[11pt]{amsart}
\usepackage{amsmath,amssymb,amsthm,array,mathtools,mathrsfs,enumitem,booktabs,float,needspace}
\usepackage[hidelinks]{hyperref}
\hypersetup{pdftitle={Incidence methods for algebraic expanders over arbitrary fields},pdfsubject={Characteristic-sensitive growth of polynomial and rational image sets over arbitrary fields},pdfkeywords={Algebraic expanders, quotients of difference sets, higher energies, point-line incidences, arbitrary fields}}
\usepackage[letterpaper,top=0.93in,bottom=0.89in,left=1.15in,right=1.15in,headheight=13pt,headsep=17pt]{geometry}
\newtheorem{theorem}{Theorem}[section]
\newtheorem{lemma}[theorem]{Lemma}
\theoremstyle{remark}
\newcommand{\F}{\mathbb F}
\newcommand{\Aff}{\operatorname{Aff}}
\newcommand{\gsim}{\mathrel{\gg}}
\newcommand{\lsim}{\mathrel{\ll}}

\newcommand{\card}[1]{\lvert #1\rvert}
\newcommand{\one}{\mathbf 1}
\newcommand{\kap}{\kappa}

\title[Incidence methods for algebraic expanders]{Incidence methods for algebraic expanders over arbitrary fields}
\author{Yao Zhi}
\email{ericyao2026@gmail.com}
\date{}
\begin{document}
\raggedbottom
\begin{abstract}
Let $A\subset F$ be finite, $N=|A|$,
$p=\operatorname{char}F$ (with $p=\infty$ in characteristic zero), and
$L=1+\log(2N)$. In the range $N\le p^{1/2}$, we establish
\[
 |AA+A|\gg N^{129/85}L^{-72/85}.
\]
This improves the long-standing $3/2$-exponent bound for $AA+A$ in positive characteristic in the range $N\le p^{1/2}$, without any assumption on the multiplicative doubling of $A$. We further establish uniform expansion bounds with exponents $129/85$ and $32/21$
for four families of three-variable product--quotient expanders, together
with refined estimates over prime fields. For $\lambda\in F^\times$ and $N\le p^{1/2}$,
our four-variable bounds include
\[
 |AA+\lambda AA|\gg N^{214/135}L^{-7/9},\qquad
 |AA-AA|\gg_\varepsilon N^{6681/4175-\varepsilon}
 \quad(\varepsilon>0).
\]
In particular, several three- and four-variable targets admit uniform
lower bounds with exponents greater than $3/2$ in the large-characteristic
range. We also prove bounds with exponent $49/38$ for shifted products and shifted quotients, and for quadratic-image
sumsets when $\operatorname{char}F\ne 2$, as well as characteristic-sensitive estimates
for $AA+\lambda A/A$
and $(A-A)^\times/(A-A)^\times$. The proofs combine higher-energy
inequalities, regularisation, and point--line incidences.
\end{abstract}
\maketitle
\clearpage
\tableofcontents

\section{Introduction}\label{sec:introduction}
\subsection{Background and earlier results on algebraic expanders}
A central theme in additive combinatorics is the growth of sets obtained
by combining the additive and multiplicative operations of a field.
Here we use the term \emph{algebraic expander problems} broadly to
include growth questions for image sets of both polynomial and rational maps,
with rational expressions evaluated only where their denominators are nonzero.
Unlike the classical sum--product problem, which considers the larger
of $|A+A|$ and $|AA|$, an algebraic expander problem seeks a lower
bound for a single image, such as $A(A+1)$, $(A+1)/A$, or
$A+A^{[2]}$. The behaviour of these images depends on the ambient
field: estimates over the reals need not extend to positive
characteristic, where the size of the input relative to the
characteristic is an essential parameter.

For shifted products, Jones and Roche--Newton \cite{JRN} proved
$|A(A+1)|\gg |A|^{24/19-o(1)}$ over $\mathbb R$.
Warren \cite{Warren} obtained the exponent $11/9$, up to logarithmic
factors, for arbitrary fields under the restriction
$|A|<p^{1/4}$ when the characteristic $p$ is positive.
For $A\subset\F_p^\times$ with $|A|<p^{1/2}$, Hart, Li and Shen
\cite{HLS} established $|(A+1)/A|\gg |A|^{110/109}$.
In a related additive problem, Pham, Vinh and de Zeeuw \cite{PVZ}
proved $|A+A^{[2]}|\gg |A|^{6/5}$ in odd characteristic
when $|A|\ll p^{5/8}$. Over the reals, Stevens and Warren
\cite{SW} developed higher-energy estimates for sums of convex or
concave images; their joint sumset inequality features the exponent
$49/38$ in its balanced form. These results motivate the study of
shifted products, quotients, and quadratic images within a common
energy--incidence framework.

Three-variable expanders provide a second group of closely related
problems. In suitable positive-characteristic ranges,
Roche--Newton, Rudnev and Shkredov \cite{RRS} obtained
$N^{3/2}$-scale growth for $AA+A$, while Aksoy Yazici, Murphy,
Rudnev and Shkredov \cite{AMRS} proved a corresponding estimate for
$A(A+A)$, where $N=|A|$. The latter work also gives general bounds
for $XY+Z$ and $X(Y+Z)$, which can be applied to targets involving
reciprocal sets. Over $\mathbb R$, Stevens and Warren \cite{SW}
obtained the exponent $129/85$ for $AA+A$, and Bloom's control
method \cite{Bloom} gives a $32/21$ exponent, up to a small loss,
for a related three-variable family. For finite sets of positive
real numbers, Roche--Newton \cite{RocheNewtonRatioPlusA} proved
a stronger $3/2+1/26$ exponent, with a logarithmic loss, for
$A/A+A$. The differing field and sign assumptions are relevant
when comparing these estimates.

For four-variable images, Roche--Newton, Rudnev and Shkredov
\cite{RRS} established positive-characteristic bounds for
$AA\pm AA$ that depend on the multiplicative doubling parameter
$K=|AA|/|A|$; in the applicable range, these include
$\min\{K^{1/2}N^{3/2},p\}$ as a lower-bound scale.
Rudnev and Stevens \cite{RS} subsequently obtained
$|AA+AA|\gg |A|^{127/80-o(1)}$ over $\mathbb R$.
The mixed target $AA+\lambda A/A$ combines product and quotient
representations and is naturally considered alongside these
four-variable problems.

A different geometric aspect of expansion arises from quotients of
differences. Writing $D=(A-A)\setminus\{0\}$, Shkredov
\cite{Shkredov} showed over $\mathbb R$ that
$|D/D|\gg |D|^{13/12}(\log(2|D|))^{-1/4}$.
Over prime fields, the Cartesian-product direction theorem of
Di Benedetto, Solymosi and White \cite{DiBenedettoSolymosiWhite}
implies $|D/D|\ge |A|^2-|A|$ when $|A|^2<p$.
These two statements concern distinct regimes, but both link
quotient expansion to incidence geometry.

The main geometric tool underlying our results is Lewko's
point--line incidence theorem \cite{Lewko}, which provides a
Szemer\'edi--Trotter-type estimate over arbitrary fields with an
explicit characteristic-dependent correction term. This theorem
plays a central role in our arguments, allowing us to obtain
expansion bounds that retain their dependence on the characteristic.
In particular, the precise coefficient of the characteristic term
available over prime fields is essential for our centred directional
second-moment argument. Independently, Miao and Xie \cite{MiaoXie}
obtained an incidence estimate of the same order in the prime-field
setting, while Yu \cite{Yu} subsequently refined these approaches
using Hasse derivatives. We recall the precise formulation of Lewko's
theorem in Subsection~\ref{sec:common}.

\subsection{Main results}
For the following statements, let $F$ be a field and $A\subset F$ be
finite, with $N:=|A|$ and $L:=1+\log(2N)$. Let
$p:=\operatorname{char}F$ in positive characteristic, and put
$p:=\infty$ in characteristic zero. The remaining conventions are
collected in Subsection~\ref{sec:common}.

\medskip\noindent\textbf{(A) Shifted products, quotients, and quadratic images.}
For $\lambda\in F^\times$, apply Theorem~\ref{a:thm:main} to $A_0=A\setminus\{0,-\lambda\}$. When $N\le p^{1/2}$,
\[
 \min\{|A(A+\lambda)|,\ |(A+\lambda)/A|\}
 \gg N^{49/38}L^{-67/38}.
\]

The same exponent holds for $|A+f(A)|$ and $|A-f(A)|$ when $f$ is a quadratic polynomial with nonzero leading coefficient and $\operatorname{char}F\ne2$. The full theorems give more general power--affine and characteristic-sensitive estimates; see Section~\ref{sec:shifted}.

\medskip\noindent\textbf{(B) Three-variable product--quotient expanders.}
Put $C_+=A^\times$ and $C_-=(A^\times)^{-1}$. For $\lambda\in F^\times$ and $\sigma\in\{+,-\}$, define
\[
 M_1=|A+\lambda C_\sigma A^\times|,
 \qquad M_2=|C_\sigma(A+\lambda A)|.
\]
These include the four targets $A+\lambda AA$, $A+\lambda A/A$, $A(A+\lambda A)$, and $(A+\lambda A)/A$. In the range $p\ge N^2$ (and in characteristic zero), the bounds include
\[
 M_1\gg N^{129/85}L^{-72/85},
 \qquad M_2\gg N^{32/21}L^{-83/63}.
\]
For $AA+A$, the first of these exponents improves on the $3/2$ bound of Roche--Newton, Rudnev and Shkredov~\cite{RRS} in the range $N\le p^{1/2}$, without a multiplicative-doubling hypothesis. The full characteristic-sensitive estimates and the prime-field improvements appear in Section~\ref{sec:three-mixed}.

\medskip\noindent\textbf{(C) Four-variable expanders.}
For $\lambda\in F^\times$ and $N\le p^{1/2}$,
\[
 |AA+\lambda AA|\gg N^{214/135}L^{-7/9},
 \qquad
 |AA-AA|\gg_\varepsilon N^{6681/4175-\varepsilon}
 \quad(\varepsilon>0).
\]
The first estimate is uniform in the multiplicative doubling ratio of $A$; the second uses Bloom's control theorem. Section~\ref{sec:affine-four} gives their full characteristic-dependent forms.

\medskip\noindent\textbf{(D) Mixed product--quotient expansion.}
For $A\subset F^\times$ and $\lambda\in F^\times$, let $M=|AA+\lambda A/A|$. Theorem~\ref{new:thm:mixed} gives three characteristic-dependent alternatives. In particular,
\[
 M\gg N^{214/135}L^{-2/9}(L+N^2/p)^{-5/9}
 \qquad(p\ge N^{8/5}),
\]
and hence
\[
 M\gg N^{214/135}L^{-7/9}
 \qquad(p\ge N^2/L).
\]
Both conclusions also hold in characteristic zero.

\medskip\noindent\textbf{(E) Symmetric difference-quotient expansion.}
Let $D=(A-A)\setminus\{0\}$, $d=|D|$, $n=|A|\ge3$, and
\[
 \mathcal H_n=1+\log(4n)+\min\{4n^2/p,2n\},
\]
with the characteristic term omitted in characteristic zero. Theorem~\ref{new:thm:symmetric-difference-quotient} gives two complementary lower bounds for $s=|D/D|$:
\[
\begin{split}
 s\gg\max\Bigg\{
 &\min\left\{\frac{d^{5/6}n^{1/2}}{\mathcal H_n^{1/4}},
             \frac{p^{1/6}d^{2/3}n^{1/3}}{\mathcal H_n^{1/6}}\right\},\\[-2pt]
 &\min\left\{\frac{d^{3/2}}{n^{3/4}\mathcal H_n^{3/8}},
             \frac{p^{1/4}d}{n^{1/2}\mathcal H_n^{1/4}}\right\}
 \Bigg\}.
\end{split}
\]
In characteristic zero, only the first entry of each minimum is retained.

\subsection{Comparison with previous results}
Table~\ref{tab:earlier-expansion} compares representative specialisations of our main theorems with previously known lower bounds applicable to the corresponding expander sets. The final
column records the \emph{common} field and size hypotheses, not the full
range of our estimates; in particular, comparisons over $\mathbb R$ or
$\F_p$ are not statements about arbitrary extension fields.
The last row uses the notation of part (E) above. The parameter
$\varepsilon>0$ is fixed, and the earlier shifted-product entry suppresses
logarithmic losses.

\begin{table}[!htbp]
\centering
\caption{Comparison of our results with previously known bounds.}
\label{tab:earlier-expansion}
\begingroup
\footnotesize
\setlength{\tabcolsep}{2.1pt}
\renewcommand{\arraystretch}{1.09}
\begin{tabular}{@{}p{.22\textwidth}p{.228\textwidth}p{.32\textwidth}p{.168\textwidth}@{}}
\toprule
\textbf{Target} & \textbf{Earlier bound} & \textbf{This paper (theorem)} & \textbf{Common range} \\
\midrule
\multicolumn{4}{@{}l}{\emph{I. Comparisons in positive characteristic}}\\
\addlinespace[1pt]
$|A(A+1)|$ & $\gg N^{11/9}$ \cite{Warren} &
$\gg N^{49/38}L^{-67/38}$\newline (Thm.~\ref{a:thm:main}) & $N<p^{1/4}$ \\
$|(A+1)/A|$ & $\gg N^{110/109}$ \cite{HLS} &
$\gg N^{49/38}L^{-67/38}$\newline (Thm.~\ref{a:thm:main}) & $F=\F_p$,\newline $N\le p^{1/2}$ \\
$|A+A^{[2]}|$ & $\gg N^{6/5}$ \cite{PVZ} &
$\gg N^{49/38}L^{-67/38}$\newline (Thm.~\ref{a:thm:quadratic}) & $p>2$,\newline $N\le p^{1/2}$ \\
$|AA+A|$ & $\gg N^{3/2}$ \cite{RRS} &
$\gg N^{129/85}L^{-72/85}$\newline (Thm.~\ref{new:thm:threefamily}) & $p>2$,\newline $N\le p^{1/2}$ \\
$|A(A+A)|$ & $\gg N^{3/2}$ \cite{AMRS} &
$\gg N^{32/21}L^{-83/63}$\newline (Thm.~\ref{new:thm:threefamily}) & $p>2$,\newline $N\le p^{1/2}$ \\
$|AA+AA|$ & $\gg N^{3/2}$ \cite{RRS} &
$\gg N^{214/135}L^{-7/9}$\newline (Thm.~\ref{b:th:main}) & $p>2$,\newline $N\le p^{1/2}$ \\
$|AA-AA|$ & $\gg N^{3/2}$ \cite{RRS} &
$\gg_\varepsilon N^{6681/4175-\varepsilon}$\newline (Thm.~\ref{b:th:diff}) & $p>2$,\newline $N\le p^{1/2}$ \\
\midrule
\multicolumn{4}{@{}l}{\emph{II. Additional comparisons and prime-field refinements}}\\
\addlinespace[1pt]
$|A+\lambda A/A|$ & $\gg N^{3/2}$\newline \cite[Cor.~20]{AMRS} &
$\gg N^{129/85}L^{-72/85}$\newline (Thm.~\ref{new:thm:threefamily}) & $p>2$,\newline $p\ge N^2$ \\
$|A+A/A|$ ($\lambda=1$) & $\gg N^{3/2+1/26}$\newline $\cdot(\log(2N))^{-5/6}$\newline \cite{RocheNewtonRatioPlusA} &
$\gg N^{129/85}L^{-72/85}$\newline (Thm.~\ref{new:thm:threefamily}) & $F=\mathbb R$,\newline $A\subset\mathbb R_{>0}$ \\
$|(A+\lambda A)/A|$ & $\gg N^{3/2}$\newline \cite[Cor.~20]{AMRS} &
$\gg N^{32/21}L^{-83/63}$\newline (Thm.~\ref{new:thm:threefamily}) & $p>2$,\newline $p\ge N^2$ \\
$M_1,M_2$ over $\F_p$ & $M_i\gg N^{3/2}$\newline \cite[Cor.~20]{AMRS} &
$M_i\gg\mathcal P_i$ ($i=1,2$)\newline (Thm.~\ref{new:thm:primefamily}) & $F=\F_p$,\newline $p\ge N^{3/2}$ \\
$|AA+\lambda A/A|$ & $\gg\min\{N^{3/2},p\}$\newline \cite[Cor.~20]{AMRS} &
$\gg N^{214/135}L^{-7/9}$\newline (Thm.~\ref{new:thm:mixed}) & $A\subset F^\times$,\newline $p>2$,\newline $p\ge N^2/L$ \\
$|D/D|$ & $\ge n^2-n$\newline \cite{DiBenedettoSolymosiWhite} ($\F_p$);\newline
$\gg d^{13/12}$\newline $\cdot(\log(2d))^{-1/4}$\newline \cite{Shkredov} ($\mathbb R$) &
$\gg\max\{\mathcal B_{\rm all},\mathcal B_{\rm sub}\}$\newline (Thm.~\ref{new:thm:symmetric-difference-quotient}) & $\F_p$: $n^2<p$;\newline $\mathbb R$: $n\ge3$ \\
\bottomrule
\end{tabular}
\endgroup
\end{table}

In the prime-field row, abbreviate the estimates
of Theorem~\ref{new:thm:primefamily} as
\begin{align*}
\mathcal P_1&=\min\bigl\{N^2/H_0,\ N^{129/85}L^{-96/85},\\[-2pt]
 &\hspace{4em}N^{219/152}p^{3/76}L^{-177/152}\bigr\},\\
\mathcal P_2&=\min\bigl\{p,\ N^{32/21}L^{-107/63},\\[-2pt]
 &\hspace{4em}N^{27/19}p^{1/19}L^{-98/57}\bigr\}.
\end{align*}
The symbols $\mathcal B_{\rm all}$ and $\mathcal B_{\rm sub}$ refer to
the two explicit minima in
Theorem~\ref{new:thm:symmetric-difference-quotient}.

The less immediate positive-characteristic entries in Panel~II follow
from \cite[Corollary~20]{AMRS}: for $p>2$ and finite $X,Y,Z\subset F$,
with $M=\max\{|X|,|Y|,|Z|\}$, both $|XY+Z|$ and $|X(Y+Z)|$ are bounded
below by
\[
 \gg\min\bigl\{(|X||Y||Z|)^{1/2},\ |X||Y||Z|/M,\ p\bigr\}.
\]
For $A\subset F^\times$, choose $(X,Y,Z)=(\lambda A,A^{-1},A)$ for
$A+\lambda A/A$ and $(X,Y,Z)=(A^{-1},A,\lambda A)$ for
$(A+\lambda A)/A$; the corresponding choices of $A^\times$ cover both
signs in $M_1,M_2$ when $0\in A$. For the mixed target, fixing
$a_0\in A$ yields
\[
 (\lambda A)A^{-1}+a_0A\subset AA+\lambda A/A,
\]
which gives the stated $\min\{N^{3/2},p\}$ benchmark.

For the last row, the prime-field direction bound in
\cite[Theorem~1]{DiBenedettoSolymosiWhite} loses precisely the two
axial directions on passage to $D/D$. Likewise, removing zero from
the full difference-set quotient of \cite[Theorem~3]{Shkredov} changes
its size by one. These comparisons are confined to their stated fields.
The Hart--Li--Shen estimate assumes $A\subset\F_p^\times$; the
$3/2$ benchmark for $A(A+A)$ uses
\cite[Corollary~3(a)]{AMRS}, rather than the weaker
\cite[Theorem~14]{RRS} bound. Finally, the real-field estimate for
$A+A/A$ has a larger exponent than the bound proved here, and the
prime-field centred minima need not exceed $N^{3/2}$ for every $N$
because of logarithmic losses. Thus the table makes no universal
claim of numerical superiority or optimality.

\subsection{Proof strategy and organisation}
Our proofs combine higher-energy inequalities with point--line incidence
estimates over arbitrary fields. Mixed-control inequalities and
representation-transfer arguments convert bounds for representation counts
into growth estimates, while weighted regularisation handles uneven
multiplicities. Throughout, we retain the characteristic-dependent term
in Lewko's incidence bound. These tools yield the shifted-product,
shifted-quotient and quadratic-image results.

For the three-variable expanders, we use incidences associated with
structured affine transformations and compare additive and multiplicative
energies; a centred directional second moment gives the prime-field
refinement. For four-variable and mixed product--quotient targets, dyadic
decompositions and collision estimates control the product and quotient
representation levels. The symmetric difference quotient is treated using
Cartesian-grid ratio counts and multiplicative Pl\"unnecke--Ruzsa
inequalities.

The remaining sections are organised by target family. Section~\ref{sec:shifted}
develops the mixed-energy and regularisation tools and applies them to
shifted images. Section~\ref{sec:affine-four} establishes affine-group
estimates and four-variable expansion bounds. Section~\ref{sec:three-mixed}
treats the three-variable and mixed product--quotient expanders, and
Section~\ref{sec:diff-ratios} deals with symmetric difference quotients.

\Needspace{11\baselineskip}
\subsection{Notation and incidence estimates}\label{sec:common}
We collect the set-theoretic and asymptotic conventions used throughout.
For finite $X,Y\subset F$, write
\[
\begin{aligned}
 F^\times&:=F\setminus\{0\},&
 X^\times&:=X\cap F^\times,\\
 X^{-1}&:=\{x^{-1}:x\in X^\times\},&
 X^{[2]}&:=\{x^2:x\in X\},\\
 X+Y&:=\{x+y:x\in X,\ y\in Y\},&
 X-Y&:=\{x-y:x\in X,\ y\in Y\},\\
 XY&:=\{xy:x\in X,\ y\in Y\},&
 X/Y&:=\{x/y:x\in X,\ y\in Y^\times\}.
\end{aligned}
\]
Variables in different factors are chosen independently, and all
quotient sets exclude zero denominators. Throughout, $1/p=0$ in
characteristic zero, and alternatives containing a positive power of
$p$ are omitted from minima in that case. Logarithms are natural.
For nonnegative quantities $U,V$, our asymptotic notation means
\[
\begin{aligned}
 U\ll V&\iff U\le CV\text{ for some constant }C>0,\\
 U\gg V&\iff V\ll U,\\
 U\asymp V&\iff (U\ll V\text{ and }V\ll U).
\end{aligned}
\]
Implied constants are independent of the field and input sets unless
a subscript indicates otherwise. If $N$ denotes the size of an input
set other than $A$, this is stated locally.

We fix the counting conventions used in the proofs. For finite subsets
$X,Y$ of an abelian group $G$, written additively, and $j\ge 2$, define
\[
\begin{aligned}
 r_{X-Y}(t)&:=\bigl|\{(x,y)\in X\times Y:x-y=t\}\bigr|,\\
 E_j(X,Y)&:=\sum_{t\in G}r_{X-Y}(t)^j,
 & E_j(X)&:=E_j(X,X).
\end{aligned}
\]
In the multiplicative group $F^\times$ the corresponding definitions
for $X,Y\subset F^\times$ are
\[
\begin{aligned}
 r_{X/Y}(t)&:=\bigl|\{(x,y)\in X\times Y:x/y=t\}\bigr|,\\
 E_j^\times(X,Y)&:=\sum_{t\in F^\times}r_{X/Y}(t)^j,
 & E_j^\times(X)&:=E_j^\times(X,X).
\end{aligned}
\]
For finite $X,Y\subset F$, the sum and product representation
functions are given by
\[
\begin{aligned}
 r_{X+Y}(t)&:=\bigl|\{(x,y)\in X\times Y:x+y=t\}\bigr|,\\
 r_{XY}(t)&:=\bigl|\{(x,y)\in X\times Y:xy=t\}\bigr|.
\end{aligned}
\]
At the second-moment level,
\[
 E_2(X,Y)=\sum_{t\in F} r_{X+Y}(t)^2,\qquad
 E_2^\times(X,Y)=\sum_{t\in F^\times} r_{XY}(t)^2
 \quad (X,Y\subset F^\times).
\]
These identities do not generally extend to $j>2$. Unless indicated
otherwise, energies refer to the relevant group operation. Dyadic
levels have the form $\{t:2^k\le r(t)<2^{k+1}\}$, with empty levels omitted.

\Needspace{12\baselineskip}
The incidence estimate used throughout the paper is the following
precise form of the arbitrary-field Szemer\'edi--Trotter theorem
\cite[Theorem~1]{Lewko}.
\begin{theorem}[Lewko's point--line incidence theorem]
\label{thm:lewko-incidence}
Let $F$ be a field and let $P\subset F^2$ be a finite set of $m$
points. Let $\mathcal L$ be a family of $n$ distinct affine lines,
and put
\[
 I(P,\mathcal L):=
 \bigl|\{(u,\ell)\in P\times\mathcal L:u\in\ell\}\bigr|.
\]
Then
\begin{equation}\label{eq:commonST}
 I(P,\mathcal L)\le 3(mn)^{2/3}+m+n+\frac{2mn}{p},
\end{equation}
where $p=\operatorname{char}F$ in positive characteristic and $1/p=0$
in characteristic zero. If $0<p<\infty$ and $mn\le p^3$, the stronger
estimate
\[
 I(P,\mathcal L)\le 3(mn)^{2/3}+m+n
\]
holds. If $F=\F_p$, the coefficient $2$ of $mn/p$ in
\eqref{eq:commonST} can be replaced by $1$ without any size
restriction.
\end{theorem}

Independently, Miao and Xie \cite{MiaoXie} proved a point--line
incidence bound of the same order over prime fields using a different
polynomial-decomposition method. Yu \cite{Yu} compared the two approaches
and obtained refined characteristic-dependent estimates using Hasse
derivatives, including bounds sharp in all parameter ranges over
algebraically closed fields. We retain Lewko's explicit constants
because the coefficient of $mn/p$ over prime fields is used in the
centred directional second-moment argument. Throughout, $p$ denotes
the characteristic rather than the order of an extension field.
Passing to subsets of size comparable to $N$ changes the logarithmic
and characteristic correction factors used below by at most absolute
constants.
\section{Shifted products, quotients and quadratic images}\label{sec:shifted}
\subsection{Expansion under power-affine and quadratic relations}
We treat power--affine relations $f(x)^s=\alpha x^r+\beta$ in the multiplicative group and quadratic images in the additive group. Both admit the incidence models required by the mixed-control method.

For the input size $N$, define the correction factors
\[
 \mathcal H=L^5+L^2N^2/p,\qquad H_0=1+N^2/p.
\]

\begin{theorem}\label{a:thm:main}
Let $X,Y\subset F^\times$ have cardinality $N\ge1$. Suppose there is a bijection $f:X\to Y$ such that
\[
 f(x)^s=\alpha x^r+\beta\qquad (x\in X),
\]
where $r,s$ are fixed positive integers and $\alpha,\beta\in F^\times$. Then
\begin{equation}\label{a:eq:main}
 \min\{\card{XY},\card{X/Y}\}\gg_{r,s}
 N^{49/38}L^{-6/19}\mathcal H^{-11/38}.
\end{equation}
A direct incidence argument also gives
\begin{equation}\label{a:eq:direct}
 \min\{\card{XY},\card{X/Y}\}\gg_{r,s}
 \min\{N^{5/4},(pN)^{1/2}\}.
\end{equation}
Moreover, when $N\le p^{2/3}$,
\begin{equation}\label{a:eq:wide}
 \min\{\card{XY},\card{X/Y}\}\gg_{r,s} N^{49/38}L^{-253/38}H_0^{-3/38}.
\end{equation}

\end{theorem}

\begin{theorem}\label{a:thm:quadratic}
Suppose $\operatorname{char} F\ne2$ and $f(x)=ax^2+bx+c$ with $a\ne0$. For finite nonempty $A\subset F$,
\[
 \card{A+f(A)}\gg N^{49/38}L^{-6/19}\mathcal H^{-11/38},
 \qquad N=\card A.
\]
The same estimate holds for $A-f(A)$. Furthermore,
\[
 \min\{\card{A+f(A)},\card{A-f(A)}\}\gg\min\{N^{5/4},(pN)^{1/2}\},
\]
and, for $N\le p^{2/3}$, the same minimum satisfies
\[
 \min\{\card{A+f(A)},\card{A-f(A)}\}\gg N^{49/38}L^{-253/38}H_0^{-3/38}.
\]

\end{theorem}

\subsection{Mixed control in abelian groups}\label{sec:mixed}
We use the mixed-control method of \cite{SW} in an arbitrary abelian group; its field-dependent inputs will be supplied by incidence estimates.

Here the ambient abelian group is written additively, even when the
application concerns $F^\times$; we use the representation and energy
notation fixed in Subsection~\ref{sec:common}. For a nonempty finite
set $X$ define
\[
 D_X=\sup_{\varnothing\ne W\text{ finite}}\frac{E_3(X,W)}{\card W^2},
 \qquad \kap(X)=D_X/\card X^2.
\]
The normalised parameter satisfies $1/\card X\le\kap(X)\le1$. In multiplicative notation, higher energies use quotient representations.

\begin{lemma}[Mixed control]\label{a:lem:mixed}
Let $X,Y$ be finite nonempty sets of size at most $N$ in an abelian group. Then
\begin{equation}\label{a:eq:mixed}
 \card{X+Y}\gg
 (\card X\card Y)^{1/2}L^{-6/19}
 (\kap(X)\kap(Y))^{-11/38}.
\end{equation}
The same bound applies to $X-Y$.
\end{lemma}
\begin{proof}
Put $x=\card X$, $y=\card Y$, $S=X+Y$, $M=\card S$, $a=E_3(X)$, and
$C_{X,Y}=\sum_h r_{X-X}(h)^2r_{Y-Y}(h)$. We first construct a subset $X_*\subset X$ of cardinality $\card{X_*}\gg x$ for which, writing $e_X=E_2(X_*)$,
\begin{equation}\label{a:eq:one}
 e_X^{3/2}\ll
 \frac{LM^{5/3}}{x^{1/2}y^3}D_X^{1/6}D_Y^{1/3}a^{1/3}C_{X,Y}.
\end{equation}
\smallskip\noindent\emph{Regularisation.} Starting with $X_0=X$, call $s\in X_i+Y$ popular if $r_{X_i+Y}(s)\ge\card{X_i}y/(cLM)$ for a sufficiently large absolute constant $c$. Retain in $X_{i+1}$ those elements for which at least half of $Y$ produces popular sums. There are at most $\card{X_i}y/(cL)$ pairs producing unpopular sums, since there are at most $M$ such sums. Each discarded element accounts for more than $y/2$ of these pairs. Thus $\card{X_{i+1}}\ge(1-2/(cL))\card{X_i}$. Since $E_2(X_0)\le N^3$, within $O(L)$ steps one has $E_2(X_{i+1})\ge E_2(X_i)/2$, while $\card{X_i}\gg x$. Set $X_*=X_i$, $X'=X_{i+1}$ and denote the popular set at this step by $P_i$.

\smallskip\noindent\emph{Weighted triples.} Consider triples $(u,u',v)\in X'\times X'\times Y$ with $u+v\in P_i$. Group them by $(h,s_1,s_2)=(u-u',u+v,u'+v)$, and write $m_\gamma$ for the number of triples in a class and $w_\gamma=r_{X'-X'}(h)$ for its weight. For each pair $(u,u')\in (X')^2$, at least $y/2$ choices of $v$ are permitted. Summing the weight over the pair gives $E_2(X')$, which is comparable to $e_X$ by the stopping condition. Thus the total weight is $\gg y e_X$. Two triples belonging to the same class differ by a displacement of the form $(t,t,-t)$. For fixed $t$, the two resulting pairs in $X'$ have difference $t$, while the corresponding pair in $Y$ has difference $-t$. Removing the popularity restriction only increases the collision count, after which these choices can be counted independently. Hence
\[
 \sum_\gamma m_\gamma^2
 \le\sum_t r_{X'-X'}(t)^2r_{Y-Y}(t)
 \le C_{X,Y}.
\]

Writing $v_0(h)=r_{X'-X'}(h)^2$, the popularity condition gives
\[
 \sum_\gamma w_\gamma^2\ll\frac{LM}{xy}
 \sum_{h}\sum_{s\in S}v_0(h)r_{X_*+Y}(s+h),
\]
where the sum over $h$ is restricted to the support of $v_0$. The last sum equals a scalar product of $\one_{X_*}*v_0$ (with reflected argument) and $r_{Y-S}$. By H\"older it is at most
$\|\one_{X_*}*v_0\|_{3/2}D_Y^{1/3}M^{2/3}$.

For every finite nonempty $D$, interpolation between $\ell^1$ and $\ell^3$ gives
\[
 \|\one_{X_*}*\one_D\|_{3/2}\le x^{1/2}D_X^{1/6}\card D^{5/6}.
\]
On the dyadic level $D_t=\{h:t\le r_{X'-X'}(h)<2t\}$, one has
\[
 t^2\card{D_t}^{5/6}\le
 \min\{e_X^{5/6}t^{1/3},a^{5/6}t^{-1/2}\}.
\]
Splitting the dyadic sum at $t=a/e_X$ and summing the two resulting geometric series gives $O(e_X^{1/2}a^{1/3})$. Thus
$\|\one_{X_*}*v_0\|_{3/2}\ll x^{1/2}D_X^{1/6}e_X^{1/2}a^{1/3}$.
Combining the estimates with $(\sum m_\gamma w_\gamma)^2\le (\sum m_\gamma^2)(\sum w_\gamma^2)$ proves \eqref{a:eq:one}.

\smallskip\noindent\emph{Combining the two estimates.} Apply \eqref{a:eq:one} in both orders. Write $b=E_3(Y)$ and $E=e_X e_Y$. Put $\Theta=C_{X,Y}C_{Y,X}/(ab)\le1$ by H\"older. Multiplication and raising to the sixth power give
\[
 E^9\ll L^{12}M^{20}(xy)^{-21}(D_XD_Y)^3(ab)^8\Theta^6.
\]
Counting representations of $X_*+Y_*$ and using Cauchy--Schwarz on difference multiplicities gives $E^{1/2}\gg (xy)^2/M$. Hence
\begin{equation}\label{a:eq:unnorm}
 (xy)^{57}\ll L^{12}M^{38}(D_XD_Y)^3(ab)^8\Theta^6.
\end{equation}
Put $\Xi=ab/[D_XD_Y(xy)^2]\le1$. Substitute $ab=\Xi D_XD_Y(xy)^2$ and $D_XD_Y=\kap(X)\kap(Y)(xy)^2$ to obtain
\[
 (xy)^{19}\ll L^{12}M^{38}(\kap(X)\kap(Y))^{11}\Xi^8\Theta^6.
\]
Since $\Xi,\Theta\le1$, this proves \eqref{a:eq:mixed}. Replace $Y$ by $-Y$ for differences.
\end{proof}
\subsection{Representation bounds from incidences}
The next two lemmas bound control parameters using representations and select large subsets suitable for the mixed-control inequality.

In the following lemmas, $G$ and $H$ are abelian groups written additively, and $f:X\to Y$ is a bijection between finite subsets. The hypotheses encode the point--line estimate from Subsection~\ref{sec:common}; the required algebraic models are verified in Subsection~\ref{new:modelsref}. Specifically, suppose there is a constant $C_f$ such that for all finite $Q,R\subset G$ and $W,C\subset H$, the number of quadruples $(q,r,w,z)$ satisfying
\[
 q+r\in X,\qquad z=f(q+r)+w,\qquad(q,r,w,z)\in Q\times R\times W\times C
\]
is bounded by
\begin{equation}\label{a:eq:admissible}
 C_f\big((uvwc)^{2/3}+(u+v)(c+w)+uvwc/p\big),
\end{equation}
where $u=\card Q$, $v=\card R$, $w=\card W$, $c=\card C$.

\begin{lemma}[Transfer of representations]\label{a:lem:rep}
Under \eqref{a:eq:admissible}, suppose that each element of $B\subset X$, $m=\card B$, has at least $\tau$ representations in $Q+R$. Then
\begin{equation}\label{a:eq:rep}
 D_{f(B)}\ll_{C_f}\frac{(1+\log(2m))u^2v^2}{\tau^3}
 +\frac{m^2uv}{p\tau}.
\end{equation}
\end{lemma}
\begin{proof}
For any $W\subset H$ of size $w\ge1$, set $r(z)=r_{f(B)+W}(z)$ and $J=\sum_z r(z)^3$. On a dyadic level $C_t=\{z:t\le r(z)<2t\}$, \eqref{a:eq:admissible} bounds $\tau t\card{C_t}$. Multiply by $t^2$ and sum. By H\"older the main term is at most $O_{C_f}((1+\log(2m))^{1/3}(uvw)^{2/3}J^{2/3})$. For the linear terms use
\[
 \sum_t t^2\card{C_t}\le E_2(f(B),-W)\le mw^2,
 \qquad\sum_t t^2\ll\min(m,w)^2\le mw.
\]
For the characteristic term use the alternative bound $E_2(f(B),-W)\le m^2w$. Thus
\[
 \tau J\ll_{C_f} (1+\log(2m))^{1/3}(uvw)^{2/3}J^{2/3}
 +m(u+v)w^2+\frac{m^2uv}{p}w^2.
\]
Either the first term dominates, giving $J\ll_{C_f}(1+\log(2m))u^2v^2w^2/\tau^3$, or one of the other terms dominates. Since $m\tau\le uv$ and $\tau\le\min(u,v)$, we have $m(u+v)\tau^2\le2u^2v^2$. This absorbs the linear contribution in the same bound. Divide by $w^2$ and take the supremum over $W$.
\end{proof}

\begin{lemma}[Control decomposition]\label{a:lem:decomp}
Under \eqref{a:eq:admissible}, if $\card X=N$, then there are subsets $U,V\subset X$, each of size at least $N/2$, such that
\begin{equation}\label{a:eq:decomp}
 \kap(U)\kap(f(V))\ll_{C_f}\frac{L^5+L^2N^2/p}{N}.
\end{equation}
In particular $U$ and $V$ need not coincide or be disjoint.
\end{lemma}
\begin{proof}
Let $T\subset X$ have size $n\ge N/2$, put $\kappa=\kap(T)$, and choose a nonempty finite witness $R_0\subset G$ of size $b$ with $E_3(T,R_0)\ge\kappa n^2b^2/2$. Since $E_2(T,R_0)\le nb^2$, differences of multiplicity less than $\kappa n/4$ contribute at most $\kappa n^2b^2/4$ to $E_3$. There is therefore a dyadic level $Q$ of size $q$ and lower multiplicity $\Delta\gg\kappa n$ such that
\[
 q\Delta^3\gg \kappa n^2b^2/L,\qquad q\Delta\le nb,\qquad \Delta\le n.
\]
In particular $b\ll L\Delta^2/(\kappa n)$.

Each of the $q\Delta$ or more pairs counted on this level yields a representation $t=q_0+r_0$ of some $t\in T$ with $(q_0,r_0)\in Q\times R_0$. The degree of any $t$ is at most $q$ or $b$, but it is not necessary to bound the witness size by $N$. Discard degrees less than $q\Delta/(2n)$ and split the remaining degrees into dyadic intervals. The actual total degree is between $q\Delta$ and $2q\Delta$, so the maximum degree is at most $2q\Delta$; the ratio of maximum degree to the cutoff is at most $4n$. Thus $O(L)$ levels suffice. Some subset $T'\subset T$, $m=\card{T'}$, has at least $\tau$ representations per element and $m\tau\gg q\Delta/L$.

Applying Lemma \ref{a:lem:rep} and normalizing, we obtain
\[
 \begin{split}
 \kap(f(T'))&\ll_{C_f}\frac{Lq^2b^2}{m^2\tau^3}+\frac{qb}{p\tau}\\
 &\ll_{C_f}\frac{L^5m}{\kappa n^2}+\frac{L^2m\Delta}{p\kappa n}\\
 &\ll_{C_f}\frac{m}{\kappa n^2}(L^5+L^2n^2/p).
 \end{split}
\]
Here the first term uses $q\Delta^3\gg\kappa n^2b^2/L$; the second uses the displayed bound for $b$ and $\Delta\le n$.

Iterate the construction, removing the nonempty subsets $T'_j$ from the successive remainders $T_j$, until their union $V$ has size at least $N/2$. Let $U$ be the remainder immediately before the last removal. Thus $\card U\ge N/2$, and $U\subset T_j$ at every step. Monotonicity of unnormalized control gives $D(T_j)\ge D(U)$ because $U\subset T_j$ and $E_3(U,W)\le E_3(T_j,W)$ for every witness $W$; in other words,
$\kap(T_j)\card{T_j}^2\ge \kap(U)\card U^2$. Therefore
\[
 \kap(f(T'_j))\ll_{C_f}\frac{\card{T'_j}\mathcal H}{\kap(U)N^2}.
\]
For disjoint sets $Z_j$ of sizes $m_j$, Minkowski in $\ell^3$ followed by H\"older gives
\[
 D_{\cup_jZ_j}^{1/3}\le\sum_jD_{Z_j}^{1/3},\qquad
 \kap(\cup_jZ_j)\le\sum_j\kap(Z_j).
\]
Indeed $(\sum_j\kap(Z_j)^{1/3}m_j^{2/3})^3\le(\sum_j\kap(Z_j))(\sum_jm_j)^2$. The images of the removed sets are disjoint because $f$ is injective. Summing the previous bounds gives \eqref{a:eq:decomp}.
\end{proof}

\subsection{Weighted regularisation and the wider characteristic range}\label{sec:weighted}
Weighted regularisation sharpens the characteristic dependence by applying \eqref{a:eq:admissible} to $f$ and $f^{-1}$.

\begin{lemma}\label{a:lem:weighted}
Let $T_0,R$ be nonempty subsets of the input group with $\card{T_0},\card R\le N$ and $T_0$ contained in the domain of $f$. There exist $T'\subset T_*\subset T_0$ with $m=\card{T_*}\gg \card{T_0}$ and $\card{T'}\gg m/L$. Writing $K=E_3(T_*,R)$ and $B=E_2(T_*,R)$, for every nonempty finite $W$ in the output group one has
\begin{equation}\label{a:eq:weighted}
 \frac{E_3(f(T'),W)}{\card W^2}\ll_{C_f}
 \frac{L^2m^3N^2}{K}+\frac{mNB}{pK}\frac{E_2(f(T'),W)}{\card W}.
\end{equation}
\end{lemma}
\begin{proof}
\smallskip\noindent\emph{Selection by weighted degree.} For a current set $T_i$ assign the weight $\sum_{r\in R}r_{T_i-R}(x-r)^2$ to each $x\in T_i$. The equality of the total weight with the third energy is exact:
\[
 \sum_{x\in T_i}\sum_{r\in R}r_{T_i-R}(x-r)^2
 =\sum_q r_{T_i-R}(q)^3=K_i.
\]
Set $\eta=c/L$ for a sufficiently small absolute $c>0$. Remove elements of weight greater than $K_i/(\eta\card{T_i})$; at most an $\eta$ fraction are removed. If the retained total weight is less than $K_i/2$, pass to the retained set, whose new third energy is at most half $K_i$: all difference multiplicities can only decrease after removing input elements, and the new cubic energy is bounded by the retained old weighted mass. As $K_i\ge1$ for nonempty sets and $K_0\le N^4$, this cannot repeat more than $O(L)$ times. Choosing $c$ small keeps $\card{T_i}\ge\card{T_0}/2$. At a stopping step set $T_*=T_i$ and retain as $T'$ the remaining elements with weight at least $K_i/(4\card{T_i})$. Their total weight is at least $K_i/4$ and each weight is at most $K_i/(\eta m)$, so $\card{T'}\ge\eta m/4$. Write $r(q)=r_{T_*-R}(q)$ and $h=\max r(q)\le m$.

\smallskip\noindent\emph{Incidences for the image.} For each pair $(x,w_0)$ with $x\in T'$ and $f(x)+w_0\in C$, count the representations $x=q+r_0$, $r_0\in R$, with weight $r(q)^2$. Each pair contributes $\gg K/m$. With $Q_t=\{q:t\le r(q)<2t\}$, apply \eqref{a:eq:admissible} with $Q_t,R,W,C$ and sum with weight $4t^2$. Since $\sum_t t^2\card{Q_t}\le B$, $\sum_{Q_t\ne\varnothing}t^2\ll h^2$, and $\sum_t t^2\card{Q_t}^{2/3}\ll L^{1/3}K^{2/3}$, we obtain, writing $\rho(z)=r_{f(T')+W}(z)$,
\begin{align*}
 \frac Km\sum_{z\in C}\rho(z)&\ll_{C_f} L^{1/3}K^{2/3}(\card R\card W\card C)^{2/3}\\
 &\quad +(B+h^2\card R)(\card W+\card C)+\frac{B\card R\card W\card C}{p}.
\end{align*}
\smallskip\noindent\emph{Third-moment estimate.} To pass from the incidence estimate to the third moment, split $\rho$ into dyadic levels $C_i$ and sum with weights $2^{2i}$. With $J=\sum_z\rho(z)^3$ and $w=\card W$, the main term is at most $O_{C_f}(L^{2/3}K^{2/3}(Nw)^{2/3}J^{2/3})$. The linear terms are at most $O_{C_f}(m^3Nw^2)$. Indeed,
\[
 B+h^2\card R\ll m^2N,\quad \sum_i2^{2i}\card{C_i}\le mw^2,
 \quad \sum_{C_i\ne\varnothing}2^{2i}\ll mw.
\] The characteristic term is at most $O_{C_f}(BNw E_2(f(T'),-W)/p)$. Divide by $K/m$; absorb the $J^{2/3}$ term using Young's inequality. The linear contribution is $O_{C_f}(m^4Nw^2/K)\le O_{C_f}(m^3N^2w^2/K)$. Replace $W$ by $-W$ and divide by $w^2$ to obtain \eqref{a:eq:weighted}.
\end{proof}

\begin{lemma}\label{a:lem:wide}
Let $f:X\to Y$ be a bijection, with $\card X=\card Y=N\le p^{2/3}$, for which \eqref{a:eq:admissible} holds in both directions. Then
\[
 \min\{\card{X+Y},\card{X-Y}\}\gg_{C_f}N^{49/38}L^{-253/38}H_0^{-3/38}
\]
when both groups are the same abelian group.
\end{lemma}
\begin{proof}
Apply Lemma \ref{a:lem:weighted} with $T_0=R=X$, writing $K_A,B_A$ for the resulting energies and $V_0=f(T')$, $n_1=\card{V_0}\gg N/L$. Apply it to $f^{-1}$ with $T_0=V_0$, $R=Y$; denote the energies by $K_B,B_B$, the output by $U$, and set $V=f(U)\subset V_0$. Here $n=\card U=\card V\gg n_1/L$. We have $a=E_3(U)\le K_A$ and $b=E_3(V)\le K_B$. Put $P_3=K_AK_B$.

For $J=E_3(V_0,Y)\ge K_B$ use \eqref{a:eq:weighted} with $W=Y$ and $E_2(V_0,Y)\le N\sqrt J$. It gives
\[
 J\ll_{C_f}\frac{L^2N^7}{K_A}+\frac{N^4B_A}{pK_A}\sqrt J.
\]
As $B_A^2\le N^2K_A$, absorption yields $P_3\ll_{C_f}L^2N^7+N^{10}/p^2\ll L^2N^7$. The last step uses $N^3\le p^2$.

For arbitrary $W$, use $E_2(f(T'),W)\le m^2\card W$ in \eqref{a:eq:weighted}. With $\rho_i=B_i/(pN)$, this gives
\[
 D_V\ll_{C_f}\frac{N^5}{K_A}(L^2+\rho_A),\qquad
 D_U\ll_{C_f}\frac{n_1^3N^2}{K_B}(L^2+\rho_B).
\]
Both $\rho_i\le N^2/p$ and $\rho_A\rho_B\le\sqrt{P_3}/p^2\ll L N^{7/2}/p^2\le L N^2/p$. Consequently
\[
 D_UD_V\ll_{C_f}\frac{L^4H_0 n_1^3N^7}{P_3}.
\]
Apply \eqref{a:eq:unnorm} with $X=U,Y=V$ to their sumset, or with $X=U,Y=-V$ to their difference set (reflection leaves the controls and own energies invariant). If $M$ denotes the corresponding full output size, then
\[
 n^{114}\ll_{C_f} L^{24}H_0^3 M^{38}n_1^9N^{21}P_3^5
 \ll_{C_f} L^{34}H_0^3M^{38}n_1^9N^{56}.
\]
Using $n\gg n_1/L$ followed by $n_1\gg N/L$ gives $M^{38}\gg_{C_f}N^{49}L^{-253}H_0^{-3}$.
\end{proof}

\subsection{Applications to shifted powers and quadratic images}\label{new:modelsref}
\begin{proof}[Proof of Theorem \ref{a:thm:main}]
Use multiplicative notation in $G=H=F^\times$. Writing $x=q_0r_0$ and $z=f(x)w$, the power-affine equation implies $z^s=\alpha q_0^r w^s r_0^r+\beta w^s$. Use the points $(r_0^r,z^s)$ and the lines
\[
 Y=\alpha q_0^r w^s X+\beta w^s.
\]
Each point has at most $rs$ parameter preimages; each line has at most $rs$ preimages because its intercept determines $w^s$ and its slope then determines $q_0^r$. These bounds on fibre multiplicity remain valid even when the characteristic divides an exponent. Thus each geometric incidence corresponds to at most a bounded number of the original quadruples. Equation \eqref{eq:commonST} proves \eqref{a:eq:admissible} with a constant depending only on $r,s$.

Apply Lemma \ref{a:lem:decomp} to the representation estimate just established. In the resulting notation $U,V\subset X$, the two sets to which the mixed-control lemma is applied are $U$ and $f(V)$, which need not arise from the same subset of $X$. Since $\card U,\card{f(V)}\gg N$, Lemma \ref{a:lem:mixed} gives, separately for $Uf(V)$ and $U/f(V)$, a lower bound $\gg_{r,s}N^{49/38}L^{-6/19}\mathcal H^{-11/38}$. This proves \eqref{a:eq:main}. The inverse relation is $x^r=\alpha^{-1}f(x)^s-\alpha^{-1}\beta$, so the same incidence model, with $r,s$ exchanged, establishes \eqref{a:eq:admissible} for $f^{-1}$. Lemma \ref{a:lem:wide} proves \eqref{a:eq:wide}.

For the direct bound let $M=\card{XY}$. Use points $((xb)^r,(af(x))^s)$ and lines $Y=\alpha a^s b^{-r}X+\beta a^s$, with $a,x\in X$, $b\in Y$. The points are contained in $\{t^r:t\in XY\}\times\{t^s:t\in XY\}$. The lines have at most $N^2$ distinct elements; each incidence arises from at most a bounded number (depending on $r,s$) of triples, by recovering $a,b$ up to root choices from the line and $x$ from the first coordinate. Consequently
\[
 N^3\ll_{r,s}M^{4/3}N^{4/3}+M^2+N^2+M^2N^2/p.
\]
This implies \eqref{a:eq:direct} for the product. For the quotient $S=X/Y$, use points $((x/b)^r,(a/f(x))^{-s})$ and lines $Y=\alpha b^ra^{-s}X+\beta a^{-s}$. The same count applies. Small bounded values of $N$ are covered by the trivial bound $M\ge N$.
\end{proof}

\begin{proof}[Proof of Theorem \ref{a:thm:quadratic}]
Restrict $A$ to a subset $X$ of size at least $N/2$ on which $f$ is injective. For $x=q+r$, $z=f(x)+w$, take the point $(r,z-ar^2-br)$ and the line
\[
 Y=2aqX+aq^2+bq+c+w.
\]
The point determines $r$ through its first coordinate; the slope of the line determines $q$ because $2a\ne0$, and then its intercept determines $w$. The point and line parameterisations are therefore injective. The incidence estimate gives \eqref{a:eq:admissible} in the additive group. Lemmas \ref{a:lem:decomp} and \ref{a:lem:mixed} yield the first estimates, using monotonicity and $\card X\gg N$. For the inverse on $f(X)$, write $q+r=f(z-w)$, and take the point $(z,q-az^2-bz-c)$ and line $Y=-2awX+aw^2-bw-r$. The parameters are recovered from the point and line, so \eqref{a:eq:admissible} holds in the inverse direction. Lemma \ref{a:lem:wide} yields the weighted estimates.

For the direct estimate put $B=f(A)$, $S=A+B$, $M=\card S$ and take points $(x,t-f(x))$ with $x,t\in S$ and lines $Y=-2ayX+ay^2-by+z$ for $y\in B,z\in A$. The identity $f(x-y)+z=t$ gives $N^2\card B\ge N^3/2$ incidences, by setting $x=u+y,t=f(u)+z$ with $u\in A$. There are at most $M^2$ points and $N^2$ lines, and the incidence recovers $(u,y,z)$. Apply \eqref{eq:commonST} as in the proof of Theorem \ref{a:thm:main}. Replace $f$ by $-f$ to obtain the direct estimate for $A-f(A)$.
\end{proof}
\section{Affine-group methods and four-variable expansion}\label{sec:affine-four}
\subsection{Affine-group incidence estimates for three-variable expanders}\label{sec:three-variable}
We record affine-group energy estimates needed in Section~\ref{sec:three-mixed}; the approach is related to \cite{RudShAff} and uses \eqref{eq:commonST}.

Throughout this subsection, $E_j^+$ and $E_j^\times$ denote additive difference energies and multiplicative quotient energies, respectively.

\subsubsection{Affine-group energy and a directional moment}
Let $\mathcal L$ be a finite collection of distinct nonhorizontal, nonvertical lines. Associate each line $y=ax+b$, where $a\ne0$, with the affine transformation $g(x)=ax+b$, and define the corresponding group energy by
\[
 E_{\Aff}(\mathcal L)=\card{\{(g_1,g_2,g_3,g_4)\in\mathcal L^4:
 g_1^{-1}g_2=g_3^{-1}g_4\}}.
\]
\begin{lemma}[Affine-group incidence bound]\label{aa:lem:affine}
For finite nonempty $S,T\subset F$ and a finite nonempty family $\mathcal L$ as above, writing $s=|S|$, one has
\begin{equation}\label{aa:eq:affine}
 I(S\times T,\mathcal L)\ll |T|^{1/2}\left(s^{2/3}|\mathcal L|^{1/3}E_{\Aff}(\mathcal L)^{1/6}
 +|\mathcal L|(1+s^2/p)^{1/2}\right).
\end{equation}
\end{lemma}
\begin{proof}
For an affine map $h$ let $i(h)=|\{x\in S:h(x)\in S\}|$, and let $r(h)=|\{(g_1,g_2)\in\mathcal L^2:g_1^{-1}g_2=h\}|$. By Cauchy--Schwarz over $T$, extending the sum of squared multiplicities to all possible images in $F$, we have
\[
 I(S\times T,\mathcal L)^2\le |T|\sum_h r(h)i(h).
\]
Here $\sum_h r(h)=|\mathcal L|^2$ and $\sum_h r(h)^2=E_{\Aff}(\mathcal L)$. Applying \eqref{eq:commonST} to $S\times S$ and the lines corresponding to maps with $i(h)\ge t$ gives
\[
 |\{h:i(h)\ge t\}|\ll s^4/t^3\qquad(c(1+s^2/p)\le t\le s),
\]
where $c$ is a sufficiently large absolute constant. Indeed, the linear and characteristic terms proportional to the number of lines are absorbed, and the remaining $s^2/t$ bound is at most $s^4/t^3$. Dyadic summation yields $\sum_{i(h)>t}i(h)^2\ll s^4/t$ for such $t$ (and also when $t>s$, with an empty sum). Thus
\[
 \sum_h r(h)i(h)\ll t|\mathcal L|^2+E_{\Aff}(\mathcal L)^{1/2}s^2t^{-1/2}.
\]
Choose $t\asymp\max\{1+s^2/p,s^{4/3}E_{\Aff}(\mathcal L)^{1/3}/|\mathcal L|^{4/3}\}$ with the constant large enough for the tail estimate. If $t>s$, the elementary bound $s|\mathcal L|^2$ is also covered by this choice. The resulting bound is
\[
 \sum_h r(h)i(h)\ll s^{4/3}|\mathcal L|^{2/3}E_{\Aff}(\mathcal L)^{1/3}
 +(1+s^2/p)|\mathcal L|^2.
\]
Taking square roots proves the lemma.
\end{proof}

For a nonempty set $B\subset F^\times$ of size $b$, let $\mathcal L_B=\{y=ax+c:a,c\in B\}$, put $R_B(t)=r_{B/B}(t)$, and for $t\ne0$ define
\[
 F_B(t)=\card{\{(c_1,c_2,c_3,c_4)\in B^4:c_2-c_1=t(c_4-c_3)\ne0\}}.
\]
\begin{lemma}[Affine-energy factorisation]\label{aa:lem:factor}
With $J_b=1+\log(2b)+b^2/p$, one has
\begin{align}
 E_{\Aff}(\mathcal L_B)&=\sum_{t\ne0} R_B(t)^2\big(b^2+F_B(t)\big),\label{aa:eq:factor}\\
 \sum_{t\ne0} F_B(t)^2&\ll b^6 J_b.\label{aa:eq:direction}
\end{align}
\end{lemma}
\begin{proof}
Write $g_{a,c}(x)=ax+c$. Since
\[
 g_{a_1,c_1}^{-1}g_{a_2,c_2}(x)=\frac{a_2}{a_1}x+\frac{c_2-c_1}{a_1},
\]
equality of two relative transformations is equivalent to $a_1/a_3=a_2/a_4=t$ and $c_2-c_1=t(c_4-c_3)$. The slope parameters contribute $R_B(t)^2$, and the translation parameters contribute $b^2+F_B(t)$, giving \eqref{aa:eq:factor}. For later use, let $r(h)=r_{B-B}(h)$, an even function of $h$. Then, for every $t\ne0$,
\[
 F_B(t)=\sum_{h\ne0}r(h)r(th),\qquad
 E_2^+(B,tB)=\sum_h r(h)r(th)=b^2+F_B(t).
\]
The second identity follows by writing $x_1+ty_1=x_2+ty_2$ as $x_1-x_2=t(y_2-y_1)$ and using $r(-h)=r(h)$. In particular, $F_B(-t)=F_B(t)$ and $F_B(t^{-1})=F_B(t)$; no assumption that $B=-B$ is needed.

For a line $\ell$ let $k(\ell)=|\ell\cap(B\times B)|$. By considering points $(c_3,c_1)$ and $(c_4,c_2)$, we obtain
\[
 F_B(t)=\sum_{\ell:\,\operatorname{slope}(\ell)=t}k(\ell)(k(\ell)-1),\qquad
 F_B(t)^2\le b^2\sum_{\substack{\ell:\,\operatorname{slope}(\ell)=t\\k(\ell)\ge2}}k(\ell)^3.
\]
The last inequality follows from $\sum_{\ell:\,\operatorname{slope}(\ell)=t}k(\ell)=b^2$, where empty lines contribute zero. On the entire grid the sum of $k(\ell)(k(\ell)-1)$ over lines is at most $b^4$, by uniqueness of the line through two distinct points. The contribution of lines with $2\le k(\ell)\le c(1+b^2/p)$ to the third moment is therefore $O(b^4(1+b^2/p))$. For larger richness $u$, \eqref{eq:commonST} gives $O(b^4/u^3)$ lines of richness at least $u$ by the same absorption as in Lemma~\ref{aa:lem:affine}. Here $u\le b$ for any line in the grid. Each dyadic level contributes $O(b^4)$, and there are $O(1+\log(2b))$ such levels. Hence $\sum_{k(\ell)\ge2}k(\ell)^3\ll b^4J_b$, which proves \eqref{aa:eq:direction}.

\end{proof}

\subsection{Four-variable expansion}\label{sec:four-variable}
The sum-of-products estimate is uniform in the multiplicative doubling ratio $|AA|/|A|$, unlike earlier finite-field estimates that depend on this ratio \cite{RRS}. The proof uses dyadic product-representation levels and a common collision budget; the difference bound additionally invokes Bloom's control theorem \cite{Bloom}.

For the four-variable estimates, introduce the characteristic correction $J=L+N^2/p$.
All sets under consideration are nonempty.
\begin{theorem}\label{b:th:main}
Let $A\subset F$ be a finite nonempty set and $\lambda\in F^\times$. Then
\begin{equation}\label{b:eq:main}
 \card{AA+\lambda AA}\gg \min\left\{
 N^{214/135}L^{-2/9}J^{-5/9},\quad
 p^{4/31}N^{42/31}L^{-6/31}J^{-15/31}\right\}.
\end{equation}
The second term is omitted in characteristic zero. 
\end{theorem}

\begin{theorem}\label{b:th:diff}
Put $\alpha=2506/4175-\varepsilon$, where $0<\varepsilon<2506/4175$. Then
\begin{equation}\label{b:eq:diff}
 \card{AA-AA}\gg_\varepsilon N^{1+\alpha}J^{-\alpha}.
\end{equation}
\end{theorem}
\subsubsection{Energy and representation estimates}
We use additive energies $E_j$ and write $D(X)=D_X$, $\kappa(X)=D(X)/|X|^2$ as in Subsection~\ref{sec:mixed}. Multiplicative product representations are kept separate.

The following lemma records the representation estimate in the form needed for product levels.
\begin{lemma}\label{b:lem:repr}
Let $Q,R\subset F^\times$ have sizes $u,v$, and let $B\subset F^\times$ have size $n$. Suppose that each $b\in B$ has at least $t$ representations either as $qr$ or as $q/r$, with $(q,r)\in Q\times R$ (use the same convention for all $b$). Then
\begin{equation}\label{b:eq:rep}
 D(B)\ll \frac{(1+\log(2n))u^2v^2}{t^3}
       +\frac{n^2uv}{pt}.
\end{equation}
\end{lemma}
\begin{proof}
In the quotient case replace $R$ by $R^{-1}$. Apply Lemma~\ref{a:lem:rep} with input group $F^\times$, output group $(F,+)$, and the identity on $B$. The representation model uses points $(q,z)$ and lines $y=rx+w_0$, which recover all parameters when $z=qr+w_0$. Thus the admissibility hypothesis \eqref{a:eq:admissible} follows from \eqref{eq:commonST}, and the claimed bound follows.
\end{proof}

\subsubsection{Two bounds for a product level}
Assume without loss that $0\notin A$ (deleting zero changes $N\ge2$ by at most a factor two, and all bounds in Theorems \ref{b:th:main} and \ref{b:th:diff} are stable under such a change). The singleton cases follow directly; we henceforth take $A\subset F^\times$.
Write $M=\card{AA+\lambda AA}$. For every dyadic integer $\tau\in[1,N]$ with nonempty
\[
 P=P_\tau=\{z:\tau\le r_{AA}(z)<2\tau\},\qquad m=\card P,
\]
put $v=m\tau/N^2$. Then
\begin{equation}\label{b:eq:mass}
 0<v\le1,\qquad \tfrac12\le\sum_\tau v_\tau\le1.
\end{equation}
Applying Lemma \ref{b:lem:repr} to $P$, with $Q=R=A$, gives
\begin{equation}\label{b:eq:DP}
 D(P)\ll L N^4/\tau^3+m^2N^2/(p\tau)\ll JN^4/\tau^3,
\end{equation}
since $m\tau\le N^2$. As $P+\lambda P\subset AA+\lambda AA$, Lemma \ref{a:lem:mixed} and $m=vN^2/\tau$ imply
\begin{equation}\label{b:eq:branch1}
 M\gg a_0\tau^{-8/19}v^{41/19},\qquad
 a_0=N^2L^{-6/19}J^{-11/19}.
\end{equation}

For the second bound, let $B\subset A$ consist of those $b$ for which at least $vN/2$ elements $a\in A$ satisfy $ab\in P$. The total number of such pairs before restriction is at least $m\tau=vN^2$, so $n=\card B\ge vN/2$. Each $b\in B$ has at least $vN/2$ representations as $z/a$ with $z\in P$, $a\in A$. By Lemma \ref{b:lem:repr},
\[
 D(B)\ll \frac{Lm^2}{v^3N}+\frac{n^2m}{pv}\ll\frac{Jm^2}{v^3N},
\]
where $n\le N$ and $m\ge vN$ (using $\tau\le N$) justify the final comparison. In particular, writing $F_3=E_3(B)$,
\begin{equation}\label{b:eq:F}
 F_3\ll Jn^2m^2/(v^3N).
\end{equation}
\begin{lemma}[Difference-level summation]\label{b:lem:second}
Let $B\subset A$ have size $n$, $F_3=E_3(B)$, and $M=\card{AA+\lambda AA}$. Then at least one of
\begin{align}
 M&\gg n^3N^{1/10}F_3^{-1/2},\label{b:eq:secondmain}\\
 M^3&\gg pn^6/F_3\label{b:eq:secondchar}
\end{align}
holds; in characteristic zero only \eqref{b:eq:secondmain} is needed.
\end{lemma}
\begin{proof}
Put $S=B+\lambda B$, $s=\card S$, $u=\card{AS}\le M$, and $E=E_2(B)$. By Cauchy--Schwarz on $B+\lambda B$ and then on difference multiplicities, $E\ge n^4/s$. If $S=\{0\}$, we have $n=1$ and the first assertion follows from $M\ge N$. Otherwise $u\ge\max(N,s)$: a nonzero element of $S$ gives a scaled copy of $A$, and a fixed nonzero element of $A$ gives a copy of $S$. Also $u\le Ns$.

\smallskip\noindent\emph{Retained difference levels.} For dyadic $\Delta$, let $D_\Delta=\{h:\Delta\le r_{B-B}(h)<2\Delta\}$, $d=\card{D_\Delta}$ and $e_\Delta=d\Delta^2$. Keep only levels with $\Delta\ge E/(8n^2)$. The discarded contribution to $E$ is at most $E/4$ by $\sum_h r_{B-B}(h)=n^2$. On retained levels $d\le n^2/\Delta\le 8s$, and $\sum e_\Delta\gg E$.

\smallskip\noindent\emph{A common collision bound.} For such a level count triples $(b_1,b_2,b_3)\in B^3$ with $b_1-b_2\in D_\Delta$. Send them to
\[
 (h,x,y)=(b_1-b_2,b_1+\lambda b_3,b_2+\lambda b_3).
\]
There are at least $nd\Delta$ triples. Let $C_\Delta$ be the sum of the squared fibre sizes for the triples on this level. Two triples with the same image differ by $(t,t,-t/\lambda)$. Since the levels $D_\Delta$ are disjoint and the first coordinate of every image is $h$, these collision counts are disjoint parts of a single collision count. Thus
\begin{equation}\label{b:eq:collisionbudget}
 \sum_\Delta C_\Delta\le \sum_t r_{B-B}(t)^2r_{B-B}(-t/\lambda)\le F_3,
\end{equation}
where the last inequality is H\"older. In particular $C_\Delta\le F_3$ on each level. Cauchy--Schwarz gives the level-specific inequality
\begin{equation}\label{b:eq:trip}
 n^2d^2\Delta^2\le C_\Delta\mathcal N,\qquad
 \mathcal N=\card{\{(h,x,y)\in D_\Delta\times S^2:x-y=h\}}.
\end{equation}
\smallskip\noindent\emph{Incidences.} The quantity $\mathcal N$ can be bounded by a point--line incidence argument. Each solution and each $a\in A$ gives an incidence between a point $(ax,h)\in AS\times D_\Delta$ and a line $Y=a^{-1}X-y$, for $a\in A,y\in S$. These lines are distinct, and an incidence recovers all the parameters. Thus
\[
 N\mathcal N\ll (udNs)^{2/3}+ud+Ns+udNs/p.
\]
The linear terms are absorbed into the first term: indeed $ud\le8Ns^2\le 8(Ns)^2$ and $Ns\le u^2\le(ud)^2$ (nonempty levels). Dividing by $N$ and using $u\le M$ gives
\begin{equation}\label{b:eq:inc2}
 \mathcal N\ll N^{-1/3}(Mds)^{2/3}+Mds/p.
\end{equation}
\smallskip\noindent\emph{Summation over levels.} Combine \eqref{b:eq:trip} and \eqref{b:eq:inc2}. Put
\[
 Q=F_3N^{-1/3}(Ms)^{2/3}/n^2.
\]
On each level at least one of $d^{4/3}\Delta^2\ll Q$ and $e_\Delta\ll C_\Delta Ms/(pn^2)$ holds. Assign overlapping levels to either class. In the first class,
\[
 e_\Delta\ll\min\{Q^{3/4}\Delta^{1/2},F_3/\Delta\}.
\]
The second bound holds on every level by the definition of $F_3$. Summing the displayed minimum over dyadic $\Delta$ on the two sides of $\Delta=F_3^{2/3}Q^{-1/2}$ gives $O(Q^{1/2}F_3^{1/3})$. For the characteristic-dominated levels, we use the common collision budget \eqref{b:eq:collisionbudget} instead of applying the bound $C_\Delta\le F_3$ separately on each level. Their total contribution is therefore $O(F_3Ms/(pn^2))$, without an additional logarithm from the number of levels. It follows that
\[
 \frac{n^4}{s}\le E\ll \frac{F_3^{5/6}(Ms)^{1/3}}{nN^{1/6}}
       +\frac{F_3 Ms}{pn^2}.
\]
Since $s\le M$ (fix a nonzero multiplier in $A$), if the first term dominates then $M\gg n^3N^{1/10}F_3^{-1/2}$. If the second dominates, then $M^3\gg pn^6/F_3$. This proves the lemma.
\end{proof}

Inserting \eqref{b:eq:F}, $n\ge vN/2$ and $m=vN^2/\tau$ in Lemma \ref{b:lem:second}, we find that on each product level at least one of
\begin{align}
 M&\gg b_0\tau v^{5/2},& b_0&=N^{3/5}J^{-1/2},\label{b:eq:branch2}\\
 M&\gg c_0\tau^{2/3}v^{5/3},& c_0&=p^{1/3}N^{1/3}J^{-1/3}\label{b:eq:branch3}
\end{align}
holds, with the latter branch absent in characteristic zero. The estimate \eqref{b:eq:branch1} remains available on every level, regardless of which alternative holds.

\subsubsection{Geometric summation over product levels}
\begin{lemma}\label{b:lem:geometric}
Let $a,b,r,s,h,j>0$. If nonnegative dyadic weights $v_t$ on any collection of dyadic numbers $t$ satisfy
\[
 M\ge c a t^{-r} v_t^h,\qquad M\ge c b t^s v_t^j
\]
for a constant $c>0$, then with $T=a^{s/(r+s)}b^{r/(r+s)}$,
\begin{equation}\label{b:eq:geom}
 \sum_t v_t\ll_{r,s,h,j,c} (M/T)^{1/h}+(M/T)^{1/j}.
\end{equation}
\end{lemma}
\begin{proof}
At $t_0=(a/b)^{1/(r+s)}$ the two coefficients coincide and equal $T$. Use the first inequality below $t_0$ and the second above it; each bound is a geometric series. The argument applies to arbitrary subcollections of dyadic numbers; in particular, the balancing point $t_0$ need not belong to the collection.
\end{proof}
\begin{proof}[Proof of Theorem \ref{b:th:main}]
Separate the product levels into two classes according to which of \eqref{b:eq:branch2}, \eqref{b:eq:branch3} is valid, assigning any overlap arbitrarily. The first inequality \eqref{b:eq:branch1} holds in both classes. For the first class, Lemma \ref{b:lem:geometric} gives the balancing scale
\[
 T_1=a_0^{19/27}b_0^{8/27}=N^{214/135}L^{-2/9}J^{-5/9}.
\]
For the second class it gives
\[
 T_2=a_0^{19/31}c_0^{12/31}=p^{4/31}N^{42/31}L^{-6/31}J^{-15/31}.
\]
More explicitly, the total mass of the two classes is at most a constant times
\[
 (M/T_1)^{19/41}+(M/T_1)^{2/5}
 +(M/T_2)^{19/41}+(M/T_2)^{3/5},
\]
with the second class omitted when $p=\infty$. By \eqref{b:eq:mass}, $M\gg\min(T_1,T_2)$.
This is exactly \eqref{b:eq:main}.
\end{proof}

\begin{proof}[Proof of Theorem \ref{b:th:diff}]
Bloom's Theorem 1 in \cite{Bloom}, formulated in an arbitrary abelian group, asserts $\card{X-X}\gg_\varepsilon \card X\kappa(X)^{-\alpha}$ for $\alpha=2506/4175-\varepsilon>0$. Apply it to each $P_\tau$. By \eqref{b:eq:DP}, $\kappa(P_\tau)\ll J/(v_\tau^2\tau)$. Consequently, with $M_-=\card{AA-AA}$,
\[
 M_-\gg_\varepsilon m (v_\tau^2\tau/J)^\alpha
     =N^2J^{-\alpha}\tau^{\alpha-1}v_\tau^{1+2\alpha}.
\]
Solving for $v_\tau$ and summing over all dyadic $1\le\tau\le N$ gives
\[
 \tfrac12\le\sum_\tau v_\tau\ll_\varepsilon
 (M_-J^\alpha/N^2)^{1/(1+2\alpha)}N^{(1-\alpha)/(1+2\alpha)}.
\]
Because $\alpha<1$, the resulting dyadic sum is geometric. This yields \eqref{b:eq:diff} without a logarithmic loss from the number of product levels. When $N\le\sqrt p$ use $J\asymp L$, then absorb the logarithm into an arbitrarily small power loss (choose $\varepsilon<\eta$).
\end{proof}

\section{Three-variable and mixed product--quotient expanders}\label{sec:three-mixed}
\subsection{Three-variable expansion and centred prime-field directions}\label{new:sec:threefamily}
We use the correction factors $\mathcal H,H_0,J$ defined earlier. Put
$C_+=A^\times$ and $C_-=(A^\times)^{-1}$.
For $\lambda\in F^\times$ and $\sigma\in\{+,-\}$, set
\[
 M_1=|A+\lambda C_\sigma A^\times|,\qquad
 M_2=|C_\sigma(A+\lambda A)|.
\]

\begin{theorem}[Uniform three-variable bounds]\label{new:thm:threefamily}
Over an arbitrary field,
\begin{align}
 M_1&\gsim \min\{N^{129/85}L^{-6/85}\mathcal H^{-11/85}J^{-11/85},\;N^2/H_0\},\label{eq:M1}\\
 M_2&\gsim N^{32/21}L^{-6/63}\mathcal H^{-11/63}J^{-22/63}.\label{eq:M2}
\end{align}
In the range $p\ge N^{3/2}$, we additionally have
\begin{align}
 M_1&\gsim\min\{N^2/H_0,\;N^{129/85}L^{-1}J^{-11/85},\;
 N^{219/152}p^{3/76}L^{-79/76}J^{-1/8}\},\label{eq:three-p1}\\
 M_2&\gsim\min\{N^{32/21}L^{-85/63}J^{-22/63},\;
 N^{27/19}p^{1/19}L^{-79/57}J^{-1/3}\}.\label{eq:three-p2}
\end{align}
In characteristic zero, branches involving positive powers of $p$ are omitted.
\end{theorem}

\begin{theorem}[Prime-field improvement]\label{new:thm:primefamily}
For $F=\F_p$ and $p\ge N^{3/2}$,
\begin{align}
 M_1&\gsim\min\{N^2/H_0,\;N^{129/85}L^{-96/85},\;
 N^{219/152}p^{3/76}L^{-177/152}\},\\
 M_2&\gsim\min\{p,\;N^{32/21}L^{-107/63},\;
 N^{27/19}p^{1/19}L^{-98/57}\}.
\end{align}
\end{theorem}

\begin{lemma}[Centred directional second moment in $\F_p$]\label{new:lem:centered}
Let $B\subset\F_p$, $b=|B|$, and $G_B(t)=E_2^+(B,tB)$. If $b^2\ge p$, then
\begin{equation}\label{eq:centered}
 \sum_{t\in\F_p^\times}\left(G_B(t)-\frac{b^4}{p}\right)^2
 \lsim b^6(1+\log(2b)).
\end{equation}
For $B\subset\F_p^\times$, write $R_B(t)=r_{B/B}(t)$ and $E_j^\times(B)=\sum_tR_B(t)^j$. Without a density restriction,
\begin{align}
 \sum_t R_B(t)G_B(t)&\lsim b^6/p+b^3(1+\log(2b))^{1/2}E_2^\times(B)^{1/2},\\
 \sum_t R_B(t)^2G_B(t)&\lsim (b^4/p)E_2^\times(B)+b^3(1+\log(2b))^{1/2}E_4^\times(B)^{1/2}.
\end{align}
\end{lemma}
\begin{proof}
Write $\mu=b^2/p$, and let $k_\ell=|\ell\cap(B\times B)|$ and
$d_\ell=k_\ell-\mu$ for each affine line in $\F_p^2$, including empty
lines. For $t\ne0$, write $\mathcal P_t$ for the $p$ lines
$x+ty=c$, $c\in\F_p$ (a parallel class of slope $-t^{-1}$).
These lines partition $B\times B$, and their intersection counts
are exactly the representation multiplicities of $B+tB$. Thus
$\sum_{\ell\in\mathcal P_t}k_\ell=b^2$ and
$\sum_{\ell\in\mathcal P_t}d_\ell=0$. Consequently
\[
 V_t:=G_B(t)-b^4/p
   =\sum_{\ell\in\mathcal P_t}d_\ell^2
   =\sum_{\ell\in\mathcal P_t}k_\ell d_\ell\ge0.
\]
Since $t\mapsto -t^{-1}$ permutes the nonzero elements of $\F_p$, the sum over $t\ne0$ is equivalent to a sum over the corresponding nonzero slopes.
Using $d_+=\max(d,0)$ and Cauchy--Schwarz gives
\begin{equation}\label{eq:centered-cs}
 V_t^2\le b^2\sum_{\ell\in\mathcal P_t}k_\ell(d_\ell)_+^2
 =b^2\sum_{\ell\in\mathcal P_t}
       \bigl(\mu(d_\ell)_+^2+(d_\ell)_+^3\bigr).
\end{equation}
Count all $p(p+1)$ affine lines, including the vertical ones. Each grid
point is on $p+1$ lines, and every ordered pair of distinct grid points
is on one line. Hence
\[
 \sum_\ell k_\ell=(p+1)b^2,\qquad
 \sum_\ell k_\ell^2=b^4+pb^2,\qquad
 \sum_\ell d_\ell^2=pb^2-b^4/p\le pb^2.
\]
For $u\ge C$ and $\mathcal L_u=\{\ell:d_\ell\ge u\}$, the prime-field
incidence bound with its exact characteristic coefficient $1$ gives
\[
  (u-1)|\mathcal L_u|
  \le3(b^2|\mathcal L_u|)^{2/3}+b^2.
\]
For $u\le b$ this implies
$|\mathcal L_u|\ll b^4/u^3+b^2/u\ll b^4/u^3$;
for $u>b$ the level is empty. A dyadic decomposition therefore gives
$\sum_\ell(d_\ell)_+^3\ll b^4(1+\log(2b))$ above the fixed cutoff.
Below it, $(d_\ell)_+^3\ll(d_\ell)_+^2$, whose total is
$O(pb^2)=O(b^4)$ when $b^2\ge p$. Summing
\eqref{eq:centered-cs} over $t\ne0$, we conclude
\[
 \sum_{t\ne0}V_t^2
 \ll b^2\left(\mu\,pb^2+b^4(1+\log(2b))\right)
 \ll b^6(1+\log(2b)).
\]
For the weighted estimates in the dense range, write
$G_B(t)=b^4/p+V_t$ and apply Cauchy--Schwarz. This gives
\[
 \sum_tR_B(t)G_B(t)\le b^6/p+E_2^\times(B)^{1/2}
       \Big(\sum_tV_t^2\Big)^{1/2},
\]
and similarly
\[
 \sum_tR_B(t)^2G_B(t)\le (b^4/p)E_2^\times(B)
       +E_4^\times(B)^{1/2}\Big(\sum_tV_t^2\Big)^{1/2}.
\]
For the sparse range, let $F_B(t)=G_B(t)-b^2$; this is the non-diagonal
quantity in Lemma~\ref{aa:lem:factor}. Since $J_b\asymp1+\log(2b)$,
that lemma gives $\sum_tF_B(t)^2\ll b^6(1+\log(2b))$. Hence
\[
\begin{aligned}
 \sum_tR_B(t)G_B(t)
 &\le b^4+E_2^\times(B)^{1/2}\Big(\sum_tF_B(t)^2\Big)^{1/2}\\
 &\ll b^3(1+\log(2b))^{1/2}E_2^\times(B)^{1/2}.
\end{aligned}
\]
Here $E_2^\times(B)\ge R_B(1)^2=b^2$. Similarly,
\[
\begin{aligned}
 \sum_tR_B(t)^2G_B(t)
 &\le b^2E_2^\times(B)
   +E_4^\times(B)^{1/2}\Big(\sum_tF_B(t)^2\Big)^{1/2}\\
 &\ll b^3(1+\log(2b))^{1/2}E_4^\times(B)^{1/2}.
\end{aligned}
\]
In this case, $E_2^\times(B)\le bE_4^\times(B)^{1/2}$ by Cauchy--Schwarz.
The two weighted estimates therefore hold in both the dense and sparse regimes.
\end{proof}

\begin{proof}[Proof of Theorems~\ref{new:thm:threefamily} and \ref{new:thm:primefamily}]
\emph{Nonzero reduction and the ambient field.} If $0\in A$, first replace $A$ by $A^\times$, retaining at least half the elements when $N\ge2$. Each of the four target sets for the reduced input is contained in the corresponding original target. The parameters $L,\mathcal H,H_0,J$ change by at most absolute factors under this replacement, and $p\ge N^{3/2}$, when assumed, still holds for the smaller cardinality. The cases of bounded cardinality are covered by adjusting the implied constants. We relabel the reduced set and its cardinality as $A,N$; henceforth all multiplicative-group operations take place in $F^\times$.

We first record the large-energy estimates. Choose $B\subset A$ with $b=|B|\asymp N$. For $M_1$, use the line family
\[
 \mathcal L_{B,\sigma}=\{y=\lambda a^\sigma x+c:a,c\in B\},
 \qquad a^+=a,\quad a^-=a^{-1}.
\]
It consists of $b^2$ distinct nonhorizontal, nonvertical lines. For every $x\in A$, evaluating any such line gives $\lambda a^\sigma x+c\in M_1$, with $a,c\in B$. Consequently the grid $A\times M_1$ supplies at least $Nb^2$ incidences. Inversion preserves the quotient representation function in the sense $R_{B^{-1}}(t)=R_B(t^{-1})=R_B(t)$, and a fixed nonzero scaling of every slope also preserves its multiplicities. Because the intercepts remain $B$, the affine-energy factorisation and directional estimate of Lemma~\ref{aa:lem:factor} therefore apply to both signs. Applying Lemma~\ref{aa:lem:affine} to that grid gives
\begin{equation}\label{eq:largeM1}
 M_1\gsim\min\{N^{13/6}E_3^\times(B)^{-1/6}J^{-1/6},\ N^2/H_0\}.
\end{equation}
For $M_2$, let $\nu(z)$ count representations $z=a^\sigma(b_1+\lambda b_2)$ with $a,b_1,b_2\in B$. Since $B\subset A$, its support is contained in the corresponding target defining $M_2$, and $\sum_z\nu(z)=b^3$. In a collision between two such representations, fixing the ratio $t=(a'/a)^\sigma$ gives $R_B(t)$ choices for $(a,a')$ (also when $\sigma=-$, by symmetry of quotient multiplicities). The remaining equation is
\[
 b_1+\lambda b_2=t(b_3+\lambda b_4),
 \quad\text{or}\quad b_1-tb_3=\lambda(tb_4-b_2).
\]
Writing $r_{X-Y}$ for the difference representation function, its number of solutions is
\begin{align*}
 &\sum_h r_{B-tB}(h)\,r_{\lambda(tB-B)}(h)\\
 &\quad\le\left(\sum_h r_{B-tB}(h)^2\right)^{1/2}
      \left(\sum_h r_{\lambda(tB-B)}(h)^2\right)^{1/2}\\
 &\quad=E_2^+(B,tB)=G_B(t).
\end{align*}
The last equality uses the fact that scaling by $\lambda\ne0$ and reflection preserve the second energy. The domination is therefore valid for every $\lambda\ne0$, without requiring $\lambda=1$, and
\[
 \sum_z\nu(z)^2\le\sum_tR_B(t)G_B(t).
\]
For arbitrary characteristic, write $G_B(t)=b^2+F_B(t)$ as justified by the difference-function identity below and Lemma~\ref{aa:lem:factor}. Cauchy--Schwarz and that
lemma give
\[
 \begin{aligned}
 \sum_tR_B(t)G_B(t)
 &\le b^4+E_2^\times(B)^{1/2}\Big(\sum_tF_B(t)^2\Big)^{1/2}\\
 &\ll b^3J^{1/2}E_2^\times(B)^{1/2}.
 \end{aligned}
\]
Since $E_2^\times(B)^2\le b^2E_3^\times(B)$, another application of
Cauchy--Schwarz on the support of $\nu$ yields
\begin{equation}\label{eq:largeM2}
 M_2\gsim N^{5/2}E_3^\times(B)^{-1/4}J^{-1/2}.
\end{equation}

For the complementary small-energy estimates, apply the control-decomposition Lemma~\ref{a:lem:decomp} to the identity between the multiplicative and additive groups on $A\subset F^\times$. It gives $U,V\subset A$, each of size comparable to $N$, with
\[
 \kappa_\times(U)\kappa_+(V)\lsim \mathcal H/N.
\]
Put $\kappa=\kappa_\times(U)$ and fix $a_0\in A$. In the four cases, the target contains respectively
\[
 V+\lambda a_0 V,\qquad V+(\lambda/a_0)V,\qquad
 a_0(V+\lambda V),\qquad a_0^{-1}(V+\lambda V).
\]
All four containments use independent variables. In the second target, the denominator is fixed; in the final two, the outer multiplier is fixed. Consequently, each displayed subset is, up to a nonzero dilation, a sum of two homothetic copies of $V$. Applying the mixed-control Lemma~\ref{a:lem:mixed} in the \emph{additive} group therefore gives, for $i=1,2$,
\begin{equation}\label{eq:smallMi}
 M_i\gsim N^{30/19}L^{-6/19}\mathcal H^{-11/19}\kappa^{11/19}.
\end{equation}
Since $E_3^\times(U)\lsim N^4\kappa$, taking $B=U$ in \eqref{eq:largeM1}--\eqref{eq:largeM2} gives
\[
 M_1\gsim\min\{N^{3/2}\kappa^{-1/6}J^{-1/6},N^2/H_0\},
 \qquad
 M_2\gsim N^{3/2}\kappa^{-1/4}J^{-1/2}.
\]
Combining the nonconstant branches with \eqref{eq:smallMi}, with weights $(66/85,19/85)$ for $M_1$ and $(44/63,19/63)$ for $M_2$ (the large-energy bound being listed first), eliminates $\kappa$ and yields \eqref{eq:M1}--\eqref{eq:M2}. The cap $N^2/H_0$ is retained for $M_1$. When $p\ge N^2$ or $p=\infty$, the bounds simplify to the powers recorded in the introduction.

For the complementary range $p\ge N^{3/2}$, use the weighted-transfer Lemma~\ref{a:lem:weighted} with the multiplicative-to-additive identity. It gives $A_2\subset A_1\subset A$, with $|A_1|\gg N$ and $n=|A_2|\gg N/L$. Put
\[
 \mathcal E=E_3^\times(A_1,A),\quad D=D_+(A_2),\quad F_3=E_3^+(A_2).
\]
The weighted-transfer estimate, first with arbitrary auxiliary set and then with $W=A_2$, gives
\[
 \mathcal E D\lsim N^5\big(L^2+\mathcal E^{1/2}/p\big),
 \qquad
 \mathcal E F_3\lsim L^2N^7+(N^5/p)(\mathcal E F_3)^{1/2}.
\]
Young's inequality and $p\ge N^{3/2}$ imply $\mathcal E F_3\lsim L^2N^7$. Substitution in the refined mixed-control inequality gives
\begin{equation}\label{eq:three-small-energy-wide}
 M_i\gsim\min\left\{
 \frac{\mathcal E^{11/19}}{N^{14/19}L^{85/19}},
 \frac{p^{3/19}\mathcal E^{1/2}}{N^{14/19}L^{79/19}}
 \right\}.
\end{equation}
Take $B=A_1\subset F^\times$ (the initial nonzero reduction guarantees that the affine line slopes and all quotient energies are defined). Because $A_1\subset A$, each quotient representation count for $(A_1,A_1)$ is bounded by the corresponding count for $(A_1,A)$. Hence $E_3^\times(B)\le E_3^\times(A_1,A)=\mathcal E$. Balance \eqref{eq:three-small-energy-wide} against \eqref{eq:largeM1} and \eqref{eq:largeM2}. For $M_1$, the two balances use weights $(66/85,19/85)$ and $(3/4,1/4)$; for $M_2$, they use $(44/63,19/63)$ and $(2/3,1/3)$, again with the large-energy estimate first. These give the two characteristic-sensitive alternatives in \eqref{eq:three-p1}--\eqref{eq:three-p2}.

Finally, suppose $F=\F_p$ and $p\ge N^{3/2}$. By the difference-function identity in the proof of Lemma~\ref{aa:lem:factor}, the affine-energy
factorisation is $E_{\Aff}(\mathcal L_B)=\sum_tR_B(t)^2G_B(t)$. The second
weighted estimate in Lemma~\ref{new:lem:centered}, together with
$E_2^\times(B)\le b(E_3^\times(B))^{1/2}$ and
$E_4^\times(B)\le bE_3^\times(B)$, gives
\[
 E_{\Aff}(\mathcal L_B)
 \ll (b^5/p+b^{7/2}L^{1/2})(E_3^\times(B))^{1/2}
 \ll b^{7/2}L^{1/2}(E_3^\times(B))^{1/2},
\]
since $b\asymp N$ and $p\ge N^{3/2}$. Applying Lemma~\ref{aa:lem:affine} directly to
$A\times M_1$ and $\mathcal L_{B,\sigma}$ gives
\[
\begin{aligned}
 I(A\times M_1,\mathcal L_{B,\sigma})
 \ll M_1^{1/2}\Big(&b^{5/4}N^{2/3}L^{1/12}
                    (E_3^\times(B))^{1/12}\\
                 &+b^2(1+N^2/p)^{1/2}\Big).
\end{aligned}
\]
Since the incidence count is at least $Nb^2$ and $b\asymp N$, this gives
\[
 M_1\gg\min\{N^{13/6}(E_3^\times(B))^{-1/6}L^{-1/6},\ N^2/H_0\}.
\]
For $M_2$, the first weighted estimate in Lemma~\ref{new:lem:centered}
gives directly
\[
 \sum_z\nu(z)^2\ll b^6/p+b^{7/2}L^{1/2}
                         (E_3^\times(B))^{1/4}.
\]
Cauchy--Schwarz on the support therefore yields
\[
 M_2\gg\min\{p,\ b^{5/2}L^{-1/2}
                       (E_3^\times(B))^{-1/4}\}.
\]
Balancing these two prime-field large-energy estimates against
\eqref{eq:three-small-energy-wide} gives the exponents in
Theorem~\ref{new:thm:primefamily}: for $M_1$ the weights are
$(66/85,19/85)$ and $(3/4,1/4)$; for $M_2$ they are
$(44/63,19/63)$ and $(2/3,1/3)$. The claimed estimates follow.
\end{proof}

\subsection{Mixed product--quotient growth}\label{new:sec:mixed}
The product and quotient representation levels are treated separately using the two-set difference-level estimate below.

\begin{lemma}[Two-set difference-level estimate]\label{new:lem:twoset}
Let $A\subset F^\times$, $|A|=N$, let $\lambda\in F^\times$, and let $B,C\subset F^\times$ be nonempty. Set $n=|B|$, $c=|C|$, $S=B+\lambda C$, $s=|S|$, $F_3=E_3(B)$, $G_3=E_3(C)$. If $M\ge\max\{N,|AS|\}$, then
\begin{align}
 E_2(B)&\lsim \frac{F_3^{2/3}G_3^{1/6}(Ms)^{1/3}}{cN^{1/6}}
 +\frac{F_3^{2/3}G_3^{1/3}Ms}{pc^2},\label{eq:twosetB}\\
 E_2(C)&\lsim \frac{G_3^{2/3}F_3^{1/6}(Ms)^{1/3}}{nN^{1/6}}
 +\frac{G_3^{2/3}F_3^{1/3}Ms}{pn^2}.\label{eq:twosetC}
\end{align}
In characteristic zero the second term of each estimate is absent.
\end{lemma}
\begin{proof}
For dyadic $\Delta$, put $D_\Delta=\{h:\Delta\le r_{B-B}(h)<2\Delta\}$, $d_\Delta=|D_\Delta|$, $e_\Delta=d_\Delta\Delta^2$. For $\Delta<E_2(B)/(8n^2)$,
\[
 \sum_{\Delta<E_2(B)/(8n^2)}d_\Delta\Delta^2
 \le\frac{E_2(B)}{8n^2}\sum_h r_{B-B}(h)
 =\frac18 E_2(B).
\]
Each dyadic level contributes at most four times its lower-endpoint mass. Therefore the discarded levels account for at most half of $E_2(B)$, and the retained levels carry a positive proportion of the total energy. Since $S=B+\lambda C$, Cauchy--Schwarz gives
$E_2(B,-\lambda C)\ge n^2c^2/s$. On the other hand,
$E_2(B,-\lambda C)\le E_2(B)^{1/2}E_2(C)^{1/2}$ and
$E_2(C)\le c^3$, so $E_2(B)\ge n^4c/s^2$. Thus every retained
level has $\Delta\ge E_2(B)/(8n^2)\gg n^2c/s^2$, and
$d_\Delta\le n^2/\Delta\ll s^2/c$.

For each retained level, the map $(b,b',c_0)\mapsto (b-b',b+\lambda c_0,b'+\lambda c_0)$ has at least $cd_\Delta\Delta$ input triples. Writing $C_\Delta$ for the number of collisions on this level, we obtain a common collision budget by first considering all difference levels. Indeed, two triples with the same image
are related by
$(b,b',c_0)\mapsto (b+t,b'+t,c_0-t/\lambda)$.
The choices of the two $B$-pairs with increment $t$ and of the $C$-pair
with increment $-t/\lambda$ are independent. Thus
\[
 \sum_{\Delta\,{\rm retained}} C_\Delta
 \le\sum_{\Delta\,{\rm all}} C_\Delta
 =\sum_t r_{B-B}(t)^2 r_{C-C}(-t/\lambda)
 \le F_3^{2/3}G_3^{1/3}=:C_0,
\]
where the last inequality is H\"older. In particular, the same collision contributes to at most one difference level.
For each triple image $(h,x,y)$ and each $a\in A$, the point
$(ax,h)\in AS\times D_\Delta$ is incident to the line
$Y=a^{-1}X-y$. These $Ns$ lines are distinct, and each incidence
recovers its triple image and the parameter $a$. Put $u=|AS|$.
If $S$ contains a nonzero element, then $u\ge\max\{N,s\}$ and
$u\le Ns$. The retained-level bound gives
$d_\Delta\ll s^2/c$, while $s\le nc\le Nc$; hence
$ud_\Delta\ll N s^3/c\ll N^2s^2=(Ns)^2$ and
$Ns\le u^2\le(ud_\Delta)^2$. These inequalities also absorb the linear terms in Lewko's incidence estimate into the main term $(ud_\Delta Ns)^{2/3}$: after cubing, the required conditions are precisely $ud_\Delta\ll(Ns)^2$ and $Ns\le(ud_\Delta)^2$. Dividing
by $N$ and using $u\le M$ yields, for the number $N_\Delta$ of
triple images,
\[
 N_\Delta\lsim N^{-1/3}(Md_\Delta s)^{2/3}+Md_\Delta s/p.
\]
If $S=\{0\}$, then $n=c=1$ and the asserted energy bounds are
trivial because the lemma assumes $M\ge N$.
Cauchy--Schwarz on fibres gives
\[
 c^2d_\Delta^2\Delta^2\lsim C_\Delta\left(
 N^{-1/3}(Md_\Delta s)^{2/3}+Md_\Delta s/p\right).
\]
When the first term dominates, combine $e_\Delta\lsim F_3/\Delta$ with $e_\Delta\lsim Q^{3/4}\Delta^{1/2}$, $Q=C_0(Ms)^{2/3}/(c^2N^{1/3})$, and split the dyadic levels at
$\Delta_*=F_3^{2/3}Q^{-1/2}$, where the two bounds coincide.
The dyadic sums on the two sides of $\Delta_*$ are geometric; together they contribute
$\ll Q^{1/2}F_3^{1/3}
=F_3^{2/3}G_3^{1/6}(Ms)^{1/3}/(cN^{1/6})$. When the characteristic term dominates, sum $e_\Delta\lsim C_\Delta Ms/(pc^2)$ and use $\sum C_\Delta\le C_0$. This proves \eqref{eq:twosetB}; symmetry proves \eqref{eq:twosetC}.
\end{proof}

\begin{lemma}[Two-parameter dyadic summation]\label{new:lem:doublegeometric}
Let $\{v_t\}_{t\in\mathcal T}$ and $\{w_u\}_{u\in\mathcal U}$ be
nonnegative weights indexed by finite collections of positive dyadic
numbers. Suppose $a,b,r,s,h,k>0$, and for every pair $(t,u)$ with
$v_tw_u>0$ one has
\begin{equation}\label{eq:doublegeomassumptions}
 M\ge a(tu)^{-r}(v_tw_u)^h,
 \qquad M\ge b(tu)^s(v_tw_u)^k.
\end{equation}
Then
\begin{equation}\label{eq:doublegeomsum}
 \Bigl(\sum_t v_t\Bigr)\Bigl(\sum_u w_u\Bigr)
 \ll_{r,s,h,k}
 \left(\frac{M}{a}\right)^{\!s/(sh+rk)}
 \left(\frac{M}{b}\right)^{\!r/(sh+rk)}.
\end{equation}
In particular, if both weight sums are bounded below by a positive
absolute constant, then
\begin{equation}\label{eq:doublegeomconclusion}
 M\gg_{r,s,h,k} a^{s/(r+s)}b^{r/(r+s)}.
\end{equation}
\end{lemma}
\begin{proof}
For $\alpha,\beta>0$, the elementary dyadic estimate
\begin{equation}\label{eq:dyadicenvelope}
 \sum_{t\in2^{\mathbb Z}}\min\{P t^\alpha,Rt^{-\beta}\}
 \ll_{\alpha,\beta}
 P^{\beta/(\alpha+\beta)}R^{\alpha/(\alpha+\beta)}
\end{equation}
holds for all $P,R>0$. Indeed, split the sum at
$t_0=(R/P)^{1/(\alpha+\beta)}$ and sum the two geometric tails.

Put $\alpha=r/h$, $\beta=s/k$ and
$\theta=\beta/(\alpha+\beta)=sh/(sh+rk)$. Define
\[
 P=\max_t v_t t^{-\alpha},\quad R=\max_t v_t t^{\beta},
 \qquad
 Q=\max_u w_u u^{-\alpha},\quad S=\max_u w_u u^{\beta}.
\]
The two pointwise inequalities in \eqref{eq:doublegeomassumptions},
followed by the factorisation of the maxima, imply
\[
 PQ\le (M/a)^{1/h},\qquad RS\le (M/b)^{1/k}.
\]
By \eqref{eq:dyadicenvelope},
$\sum_t v_t\ll P^\theta R^{1-\theta}$ and
$\sum_u w_u\ll Q^\theta S^{1-\theta}$.
Multiplication gives
\[
 \Bigl(\sum_t v_t\Bigr)\Bigl(\sum_u w_u\Bigr)
 \ll(PQ)^\theta(RS)^{1-\theta}
 \le (M/a)^{\theta/h}(M/b)^{(1-\theta)/k},
\]
which is \eqref{eq:doublegeomsum}. Its inversion proves
\eqref{eq:doublegeomconclusion}.
\end{proof}

\begin{theorem}[Mixed expansion: characteristic-dependent bounds and a multiscale refinement]\label{new:thm:mixed}
Let $A\subset F^\times$ have $N\ge2$ elements, and let
$\lambda\in F^\times$. Set $M=|AA+\lambda A/A|$.
Then in positive characteristic
\begin{equation}\label{eq:mixedthree}
 M\gsim\min\{T_1,T_2,T_3\},
\end{equation}
where
\begin{align*}
 T_1&=N^{214/135}L^{-67/27}J^{-5/9},\\
 T_2&=p^{4/31}N^{42/31}L^{-67/31}J^{-15/31},\\
 T_3&=p^{12/151}N^{214/151}L^{-345/151}J^{-77/151}.
\end{align*}
In characteristic zero, or in positive characteristic when $p\ge N^{8/5}$,
one also has the logarithmically sharper characteristic-sensitive bound
\begin{equation}\label{eq:mixedintermediate}
 M\gg N^{214/135}L^{-2/9}J^{-5/9}.
\end{equation}
Furthermore, in characteristic zero or whenever $p\ge N^2/L$,
\begin{equation}\label{eq:mixedhigh}
 M\gg N^{214/135}L^{-7/9}.
\end{equation}
All implied constants are absolute.
\end{theorem}
\begin{proof}
Fix $a_0\in A$. The set $a_0A+\lambda A/a_0$ is contained in
$AA+\lambda A/A$ and has cardinality at least $N$, so $M\ge N$.
Select dyadic representation levels $P_\tau\subset AA$ and
$Q_\sigma\subset A/A$ with normalised masses
$v=|P_\tau|\tau/N^2\gg L^{-1}$ and
$w=|Q_\sigma|\sigma/N^2\gg L^{-1}$.
Such levels exist since each family has $O(L)$ nonempty levels and
$\sum_z r_{AA}(z)=\sum_z r_{A/A}(z)=N^2$; replacing the
representation counts by their dyadic lower endpoints loses at most
a factor of two. Write $m=|P_\tau|=vN^2/\tau$ and
$\ell=|Q_\sigma|=wN^2/\sigma$.

We spell out the selection of the auxiliary sets. For $b\in A$ put
\[
 d_P(b)=\bigl|\{a\in A:ab\in P_\tau\}\bigr|.
\]
Then $0\le d_P(b)\le N$ and
$\sum_{b\in A}d_P(b)\ge m\tau=vN^2$. Hence the set
\[
 B=\{b\in A:d_P(b)\ge vN/2\},\qquad n=|B|,
\]
satisfies $n\ge vN/2$. Each $b\in B$ has at least $vN/2$
representations $b=z/a$ with $(z,a)\in P_\tau\times A$.
Likewise, for $c\in A^{-1}$ put
\[
\begin{aligned}
 d_Q(c)&=\bigl|\{a\in A:ac\in Q_\sigma\}\bigr|,\\
 C&=\{c\in A^{-1}:d_Q(c)\ge wN/2\},
 \qquad c_0=|C|.
\end{aligned}
\]
The identity $\sum_{c\in A^{-1}}d_Q(c)\ge\ell\sigma=wN^2$
and the bound $d_Q(c)\le N$ give $c_0\ge wN/2$.
Each $c\in C$ consequently has at least $wN/2$
representations $c=z/a$ with
$(z,a)\in Q_\sigma\times A$. All these representations are
multiplicative quotients of nonzero elements, so
Lemma~\ref{b:lem:repr} applies separately to $B$ and $C$.
Using $m\tau\le N^2$, $\ell\sigma\le N^2$,
$m\ge vN$, $\ell\ge wN$, and $n,c_0\le N$, it yields
\begin{equation}\label{eq:mixed-moments}
\begin{aligned}
 F_3=E_3(B)&\ll \frac{Jn^2m^2}{v^3N},\\
 G_3=E_3(C)&\ll \frac{Jc_0^2\ell^2}{w^3N}.
\end{aligned}
\end{equation}
Here $J=L+N^2/p$ also absorbs the characteristic contributions
in \eqref{b:eq:rep}. In the estimates below we write $c=c_0$
for the cardinality of $C$.
Mixed control applied directly to $P_\tau+\lambda Q_\sigma$ gives
\begin{equation}\label{eq:mixone}
 M\gsim N^2L^{-6/19}J^{-11/19}(vw)^{41/38}(\tau\sigma)^{-4/19}.
\end{equation}
On the other hand, $AS\subset AA+\lambda A/A$ with $S=B+\lambda C$. Put
\[
 E_B=E_2(B),\quad E_C=E_2(C),\quad
 E_{BC}=\sum_x r_{B+\lambda C}(x)^2.
\]
Cauchy--Schwarz gives $n^2c^2/s\le E_{BC}\le(E_BE_C)^{1/2}$. By Lemma~\ref{new:lem:twoset}, write $E_B\ll X_B+Y_B$ and $E_C\ll X_C+Y_C$, where
\begin{align*}
 X_B&=F_3^{2/3}G_3^{1/6}(Ms)^{1/3}/(cN^{1/6}),&
 Y_B&=F_3^{2/3}G_3^{1/3}Ms/(pc^2),\\
 X_C&=G_3^{2/3}F_3^{1/6}(Ms)^{1/3}/(nN^{1/6}),&
 Y_C&=G_3^{2/3}F_3^{1/3}Ms/(pn^2).
\end{align*}
Thus $(E_BE_C)^{1/2}\ll (X_BX_C)^{1/2}+(X_BY_C)^{1/2}+(Y_BX_C)^{1/2}+(Y_BY_C)^{1/2}$. Since $s\le M$, the four terms yield
\begin{align*}
 M&\gsim (nc)^{3/2}N^{1/10}(F_3G_3)^{-1/4},\\
 M&\gsim p^{1/3}(nc)(F_3G_3)^{-1/6},\\
 M&\gsim p^{3/14}n^{9/7}c^{15/14}N^{1/28}F_3^{-3/14}G_3^{-5/28},\\
 M&\gsim p^{3/14}n^{15/14}c^{9/7}N^{1/28}F_3^{-5/28}G_3^{-3/14}.
\end{align*}
The representation-level bounds provide only lower estimates for $n$ and $c$, rather than two-sided comparability. From $m=vN^2/\tau$ and $\ell=wN^2/\sigma$ we obtain
\begin{equation}\label{eq:mix-size-safe}
 F_3\ll \frac{Jn^2N^3}{v\tau^2},\qquad
 G_3\ll \frac{Jc^2N^3}{w\sigma^2},\qquad
 n\gg vN,\quad c\gg wN.
\end{equation}
After inserting these upper bounds into each of the four inequalities
above, the remaining powers of $n$ and $c$ are nonnegative
(for the first two branches they are $n^1c^1$ and
$n^{2/3}c^{2/3}$, respectively; in the two mixed branches they are
$n^{6/7}c^{5/7}$ and $n^{5/7}c^{6/7}$). Accordingly, the lower bounds $n\gg vN$ and $c\gg wN$ suffice. No upper bounds of the form $n\ll vN$ or $c\ll wN$ are assumed. The resulting
four alternatives (at least one of which holds) are
\begin{align}
 M&\gsim N^{3/5}J^{-1/2}(vw)^{5/4}(\tau\sigma)^{1/2},\label{eq:mixalt1}\\
 M&\gsim p^{1/3}N^{1/3}J^{-1/3}(vw)^{5/6}(\tau\sigma)^{1/3},\label{eq:mixalt2}\\
 M&\gsim p^{3/14}N^{3/7}J^{-11/28}v^{15/14}w^{25/28}\tau^{3/7}\sigma^{5/14},\label{eq:mixalt3}\\
 M&\gsim p^{3/14}N^{3/7}J^{-11/28}v^{25/28}w^{15/14}\tau^{5/14}\sigma^{3/7}.\label{eq:mixalt4}
\end{align}
Combine \eqref{eq:mixone} with the first alternative using weights $19/27$ and $8/27$. This gives
\[
 M\gsim N^{214/135}L^{-2/9}J^{-5/9}(vw)^{61/54}\gsim T_1.
\]
Combining with the second alternative using weights $19/31$ and $12/31$ gives
\[
 M\gsim p^{4/31}N^{42/31}L^{-6/31}J^{-15/31}(vw)^{61/62}\gsim T_2.
\]
Finally, combine with the third alternative using weights $95/151$ and $56/151$. The resulting bound is
\[
 M\gsim p^{12/151}N^{214/151}L^{-30/151}J^{-77/151}
 v^{325/302}w^{305/302}\tau^{4/151}\gsim T_3,
\]
where
\[
 T_3=p^{12/151}N^{214/151}L^{-345/151}J^{-77/151}.
\]
The symmetric fourth alternative gives the same conclusion with $v,w$
interchanged. Here we used $v,w\gsim L^{-1}$, and the residual power
of $\tau$ is nonnegative and can be discarded since $\tau\ge1$.
At least one of the four collision terms dominates, so one of
$T_1,T_2,T_3$ is a lower bound for $M$. This proves \eqref{eq:mixedthree}.
The stronger bounds in the stated characteristic ranges follow from
Lemma~\ref{new:lem:branchdominance} and the multiscale argument below.
\end{proof}

\begin{lemma}[Uniform domination of the characteristic collision branches]
\label{new:lem:branchdominance}
Let $p<\infty$, $N\ge2$, and $J=L+N^2/p$. Let
$0<v,w\le1$ and $1\le\tau,\sigma\le N$.
Write $D_1,D_2,D_3,D_4$ for the expressions on the right of
\eqref{eq:mixalt1}--\eqref{eq:mixalt4}, respectively, without their
implicit constants, and set
\[
 \Theta=\frac{p^2J}{N^{18/5}}.
\]
Then
\begin{align}
 \frac{D_2}{D_1}
 &=\Theta^{1/6}\left(\frac{N^2}{\tau\sigma}\right)^{1/6}
     (vw)^{-5/12},\label{eq:dominance2}\\
 \frac{D_3}{D_1}
 &=\Theta^{3/28}\left(\frac N\tau\right)^{1/14}
     \left(\frac N\sigma\right)^{1/7}
     v^{-5/28}w^{-5/14},\label{eq:dominance3}\\
 \frac{D_4}{D_1}
 &=\Theta^{3/28}\left(\frac N\tau\right)^{1/7}
     \left(\frac N\sigma\right)^{1/14}
     v^{-5/14}w^{-5/28}.\label{eq:dominance4}
\end{align}
In particular, if $\Theta\ge1$, each $D_i$ is at least $D_1$.
More generally, if $\Theta$ has an absolute positive lower bound,
all four alternatives are at least a fixed multiple of $D_1$.
\end{lemma}
\begin{proof}
Divide the explicit right-hand sides of
\eqref{eq:mixalt2}--\eqref{eq:mixalt4} by that of
\eqref{eq:mixalt1} and collect the powers of $N,p,J,v,w,\tau,\sigma$.
This gives \eqref{eq:dominance2}--\eqref{eq:dominance4} exactly.
Every factor after the displayed power of $\Theta$ is at least one,
because $0<v,w\le1$ and $1\le\tau,\sigma\le N$.
The two conclusions follow.
\end{proof}

\begin{proof}[Proof of the refinements \eqref{eq:mixedintermediate}
 and \eqref{eq:mixedhigh}]
Use the complete dyadic families $P_\tau\subset AA$ and
$Q_\sigma\subset A/A$ of the preceding proof, with
$v_\tau=|P_\tau|\tau/N^2$ and
$w_\sigma=|Q_\sigma|\sigma/N^2$. Their normalised masses satisfy
\[
 \tfrac12\le\sum_\tau v_\tau\le1,
 \qquad \tfrac12\le\sum_\sigma w_\sigma\le1.
\]
Since each nonzero product or quotient has at most $N$
representations, all the relevant multiplicity levels satisfy
$1\le\tau,\sigma\le N$. Also $0<v_\tau,w_\sigma\le1$ on nonempty
levels. In contrast with the initial proof, we retain \emph{every}
nonempty pair of levels. From \eqref{eq:mixone} we have
\begin{equation}\label{eq:doublemixedfirst}
 M\gg a_0(\tau\sigma)^{-4/19}(v_\tau w_\sigma)^{41/38},
 \qquad a_0=N^2L^{-6/19}J^{-11/19},
\end{equation}
for each pair.

The four alternatives \eqref{eq:mixalt1}--\eqref{eq:mixalt4}
show that, for each pair, $M$ is bounded below (up to absolute
constants) by at least one of $D_1,D_2,D_3,D_4$ as defined in
Lemma~\ref{new:lem:branchdominance}.
In characteristic zero the characteristic terms are absent, and
only $D_1$ is required. In positive characteristic put
\[
 \Theta=\frac{p^2J}{N^{18/5}}.
\]
If $p\ge N^{8/5}$, then $J\ge N^2/p$, and hence
\[
 \Theta\ge\frac{pN^2}{N^{18/5}}\ge1.
\]
Lemma~\ref{new:lem:branchdominance} therefore implies that
$D_2,D_3,D_4\ge D_1$ for every pair, so the same pairwise lower
bound is valid regardless of which of the collision terms dominates:
\begin{equation}\label{eq:doublemixedsecond}
 M\gg b_0(\tau\sigma)^{1/2}(v_\tau w_\sigma)^{5/4},
 \qquad b_0=N^{3/5}J^{-1/2}.
\end{equation}
Apply Lemma~\ref{new:lem:doublegeometric} to
\eqref{eq:doublemixedfirst} and \eqref{eq:doublemixedsecond}.
With parameters $(r,s,h,k)=(4/19,1/2,41/38,5/4)$,
the lower bounds on both mass sums give
\[
 M\gg a_0^{19/27}b_0^{8/27}
  =N^{214/135}L^{-2/9}J^{-5/9}.
\]
This proves \eqref{eq:mixedintermediate} in positive characteristic
and its characteristic-zero counterpart.

To obtain the larger range for the pure logarithmic bound, suppose
$p\ge N^2/L$. We then have $J\le2L$ and, moreover,
\[
 \Theta\ge\frac{p^2L}{N^{18/5}}
 \ge\frac{N^{2/5}}{L}\gg1,
\]
where the last inequality follows from the elementary estimate
$1+\log(2N)\ll N^{2/5}$ for $N\ge2$.
Thus, by Lemma~\ref{new:lem:branchdominance}, every characteristic
alternative still bounds $D_1$ from below by an absolute multiple,
even if $p<N^{8/5}$. The two-parameter summation above remains
valid with absolute constants. Since $J\asymp L$, it gives
$M\gg N^{214/135}L^{-7/9}$, proving \eqref{eq:mixedhigh}.
In characteristic zero $J=L$, so the same conclusion follows
from the first part of the proof.
\end{proof}

\section{Symmetric difference-quotient expansion}\label{sec:diff-ratios}

\subsection{Point--line ratio preliminaries}\label{sec:ratio-preliminaries}
We use the factor $\mathcal H_n$ from item (E) of the introductory
results, with $n$ denoting the relevant cardinality parameter.

For a finite set $C\subset F$ and an affine line $\ell$, put
$r_\ell=|\ell\cap(C\times C)|$.

\begin{lemma}[Cartesian collinear triples]\label{new:lem:cartesian-triples}
If $C\subset F$ has at most $2n$ elements, then
\[
 \sum_{\ell:r_\ell\ge2}r_\ell^3\ll n^4\mathcal H_n.
\]
\end{lemma}
\begin{proof}
Put $m=|C|\le2n$. Every ordered collinear triple with its first two
points distinct determines a line. Lewko's incidence estimate
\eqref{eq:commonST} implies that, for
$t\ge C_0(1+m^2/p)$,
\[
 \#\{\ell:r_\ell\ge t\}\ll \frac{m^4}{t^3},
\]
with the characteristic term deleted when $p=\infty$. On the other hand,
uniqueness of the line through two distinct grid points gives
\[
 \sum_{\ell:r_\ell\ge2}r_\ell(r_\ell-1)\le m^4,
 \qquad \sum_{\ell:r_\ell\ge2}r_\ell^2\ll m^4.
\]
The contribution below the incidence cutoff is therefore
$O(m^4\min\{1+m^2/p,m\})$. On the remaining dyadic levels
$t\le r_\ell<2t$, the incidence tail and the pair count give a contribution
$O(m^4+m^2t^2)$. Summing over dyadic $t\le m$ proves
\[
 \sum_{\ell:r_\ell\ge2}r_\ell^3
 \ll m^4\bigl(1+\log(2m)+\min\{m^2/p,m\}\bigr)
 \ll n^4\mathcal H_n.
\]
\end{proof}

\begin{lemma}[Shifted intersection from incidences]\label{new:lem:shifted-intersection}
Let $U,V\subset F$ and $C_1,C_2\subset F^\times$ be finite and nonempty,
and let $a\ne0$. Then
\begin{align*}
 |U\cap(V+a)|\ll{}&
 \frac{(|UC_1|\,|VC_2|)^{2/3}}{(|C_1|\,|C_2|)^{1/3}}
 +\frac{|UC_1|}{|C_1|}+\frac{|VC_2|}{|C_2|}
 +\frac{|UC_1|\,|VC_2|}{p}.
\end{align*}
The last term is omitted in characteristic zero.
\end{lemma}
\begin{proof}
Use the point set $C_1^{-1}\times VC_2$ and the lines
$sx-ty=a$ indexed by $(s,t)\in UC_1\times C_2^{-1}$. Distinct parameter
pairs give distinct lines because $a\ne0$. Each representation $u=v+a$
with $u\in U$ and $v\in V$ yields $|C_1||C_2|$ incidences: for each
$c_1\in C_1$, $c_2\in C_2$, take
\[
 (x,y)=(c_1^{-1},vc_2),\qquad (s,t)=(uc_1,c_2^{-1}).
\]
Applying \eqref{eq:commonST} and dividing by $|C_1||C_2|$ gives the
claim.
\end{proof}

For $A\subset F$ with $|A|=n$, put
$D=(A-A)\setminus\{0\}$. If $X=-X\subset D$, define
\[
 \sigma_X=\sum_{x\in X}r_{A-A}(x),\qquad
 R_X=\left\{\frac{a-b}{b_1-b}:a,b,b_1\in A,\ b_1-b\in X\right\}.
\]

\begin{lemma}[Weighted ratios]\label{new:lem:weighted-ratios}
For every nonempty symmetric $X\subset D$,
\[
 R_X=1-R_X,\qquad |R_X|\gg\frac{\sigma_X^2}{n^2\mathcal H_n}.
\]
\end{lemma}
\begin{proof}
The substitution $(b,b_1)\mapsto(b_1,b)$ preserves the condition
$b_1-b\in X$ and sends a ratio $\lambda$ to $1-\lambda$, proving
$R_X=1-R_X$. Let $F_X(\lambda)$ count the triples $(a,b,b_1)$
contributing to $R_X$ with ratio $\lambda$. Then
\[
 \sum_\lambda F_X(\lambda)=n\sigma_X.
\]
Equality of two ratios implies collinearity of the corresponding three
pairs
\[
 (b,b'),\qquad (b_1,b_1'),\qquad (a,a')
\]
in $A\times A$. The first two points are distinct because the
denominators are nonzero. The correspondence is injective in the six
coordinates, so Lemma~\ref{new:lem:cartesian-triples} gives
$\sum_\lambda F_X(\lambda)^2\ll n^4\mathcal H_n$. Cauchy--Schwarz
now yields
\[
 n^2\sigma_X^2\le |R_X|\sum_\lambda F_X(\lambda)^2
 \ll |R_X|n^4\mathcal H_n,
\]
which proves the result.
\end{proof}

\subsection{Expansion of symmetric difference quotients}\label{sec:symmetric-difference-quotients}
Let $A\subset F$ have $n\ge3$ elements and put
\[
 \begin{aligned}
 D&=(A-A)\setminus\{0\},\qquad d=|D|,\qquad s=|D/D|,\\
 \mathcal H_n&=1+\log(4n)+\min\{4n^2/p,2n\}.
 \end{aligned}
\]
Then $D=-D$ and $n-1\le d\le n(n-1)$.

\begin{theorem}[Complementary symmetric quotient bounds]
\label{new:thm:symmetric-difference-quotient}
For every field $F$ and every finite $A\subset F$ with $n\ge3$,
\[
 s\gg\max\{\mathcal B_{\rm all},\mathcal B_{\rm sub}\},
\]
where
\[
 \begin{aligned}
 \mathcal B_{\rm all}&=\min\left\{
   \frac{d^{5/6}n^{1/2}}{\mathcal H_n^{1/4}},
   \frac{p^{1/6}d^{2/3}n^{1/3}}{\mathcal H_n^{1/6}}\right\},\\
 \mathcal B_{\rm sub}&=\min\left\{
   \frac{d^{3/2}}{n^{3/4}\mathcal H_n^{3/8}},
   \frac{p^{1/4}d}{n^{1/2}\mathcal H_n^{1/4}}\right\}.
 \end{aligned}
\]
In characteristic zero, retain only the first entry of each minimum.
The characteristic correction depends on $p$, not on the cardinality of
an extension field.
\end{theorem}

\subsubsection{Proof of the complementary bounds}
We record the two multiplicative Ruzsa--Pl\"unnecke estimates
used in the argument (see \cite{Ruzsa,Shkredov}).
For a nonempty finite set $D\subset F^\times$, write
$|D|=d$, $|D/D|=s$ and $K=s/d\ge1$. The full-set
threefold sum--difference inequality, applied in the abelian
group $F^\times$, gives
\[
 |DD/D|\ll K^3d=\frac{s^3}{d^2}.
\]
The large-subset form of the Pl\"unnecke--Ruzsa inequality,
with the two extra factors both equal to $D^{-1}$, gives
$X_0\subset D$ such that
\[
 |X_0|\ge\frac d2,\qquad
 |X_0D^{-1}D^{-1}|\ll K^2d=\frac{s^2}{d}.
\]
The implied constants are absolute; neither inequality imposes
a restriction on the characteristic. In the rest of the proof
$D^{-2}$ abbreviates $D^{-1}D^{-1}$.

\begin{proof}[Proof of Theorem~\ref{new:thm:symmetric-difference-quotient}]
First take $X=D$ in Lemma~\ref{new:lem:weighted-ratios}, and write
$R=R_D$. Since $\sigma_D=n(n-1)$,
\[
 |R|\gg\frac{n^2}{\mathcal H_n}.
\]
Moreover, $R\setminus\{0\}\subset D/D$. With $Z=|RD|$, we have
\[
 d\le Z\le1+|D^2/D|\ll\frac{s^3}{d^2},
\]
where we used $s\ge d$. Since $R=1-R$, Lemma~\ref{new:lem:shifted-intersection}
with $U=-R$, $V=R$, $a=-1$, and $C_1=C_2=D$ gives, using $Z\ge d$,
\[
 \frac{n^2}{\mathcal H_n}\ll |R|
 \ll \frac{Z^{4/3}}{d^{2/3}}+\frac{Z^2}{p}
 \ll \frac{s^4}{d^{10/3}}+\frac{s^6}{pd^4}.
\]
The two alternatives yield the two entries of $\mathcal B_{\rm all}$.

For the complementary bound choose $X_0$ as above and set
$X=X_0\cup(-X_0)$. Since $D=-D$, this is a symmetric subset of $D$ and
\[
 |X|\ge d/2,\qquad \sigma_X\ge|X|\ge d/2,
 \qquad |XD^{-2}|\ll s^2/d.
\]
The ratio set $R=R_X$ therefore satisfies
\[
 |R|\gg\frac{d^2}{n^2\mathcal H_n},\qquad R=1-R.
\]
As $R\setminus\{0\}\subset D/X$, the choice $Z=|RD|$ gives
\[
 d\le Z\le1+|D^2/X|=1+|XD^{-2}|\ll s^2/d.
\]
A second application of Lemma~\ref{new:lem:shifted-intersection} yields
\[
 \frac{d^2}{n^2\mathcal H_n}\ll |R|
 \ll\frac{Z^{4/3}}{d^{2/3}}+\frac{Z^2}{p}
 \ll\frac{s^{8/3}}{d^2}+\frac{s^4}{pd^2}.
\]
Solving the two alternatives gives $s\gg\mathcal B_{\rm sub}$.
Combining the full-set and large-subset bounds proves the theorem.
\end{proof}

\section*{Declaration of AI Assistance}
The author used ChatGPT, with GPT-6.0 available through a ChatGPT Plus subscription, as an auxiliary research and writing tool during parts of the exploratory development of this work. The tool was used for brainstorming possible proof strategies, checking intermediate calculations and logical dependencies, organizing arguments, and assisting with the exposition and preparation of the manuscript.

All mathematical statements, proofs, references, and final formulations appearing in this paper were independently reviewed and selected by the author. The author assumes full responsibility for the correctness, originality, and presentation of the results.

\end{document}